\documentclass[a4paper, 10pt]{amsart}
\usepackage{amsmath,latexsym}
\usepackage{amssymb}
\usepackage{color}
\usepackage{setspace}
\usepackage[latin1]{inputenc}
\usepackage{rotating}
\usepackage{fancybox}
\usepackage{eurosym}
\usepackage{anysize} 
\usepackage{graphicx}
\usepackage[all,2cell]{xy}
\usepackage{float}

\usepackage[ruled,vlined]{algorithm2e}

\newtheorem{defi}{Definition}[section]
\newtheorem{ej}{Example}[section]

\newtheorem{remark}{Remark}[section]
\newtheorem{corollary}{Corollary}[section]

\newtheorem{theo}{Theorem}[section]

\newtheorem{lemma}{Lemma}[section]

\newcommand{\Q}{\mathbb{Q}}
\newcommand{\Z}{\mathbb{Z}}
\newcommand{\N}{\mathbb{N}}
\newcommand{\R}{\mathbb{R}}

\newcommand{\grob}{\mathcal{G}}

\newcommand{\dsum}{\displaystyle\sum}
\newcommand{\dmax}{\displaystyle\max}

\newcommand{\dbigcup}{\displaystyle\bigcup}
\newcommand{\dbigcap}{\displaystyle\bigcap}
\marginsize{2.cm}{2.cm}{2.cm}{2.cm} 

\input{pst-plot}

\title{Partial Gr\"obner bases for multiobjective integer linear optimization}
\author{V\'ictor Blanco \and Justo Puerto}
\date{May 28, 2008}



\address{Departamento de Estad\'istica e Investigaci\'on Operativa, Universidad
de Sevilla, 41012 Sevilla, Spain}

\email{vblanco@us.es\\
puerto@us.es}

\thanks{This research was partially supported by Ministerio de
Educaci\'on y Ciencia under grant MTM2007-67433-C02-01
.}

\begin{document}
 \maketitle

 \begin{abstract}
In this paper we present a new methodology for solving
multiobjective integer linear programs using tools from algebraic
geometry.
 We introduce the concept of partial Gr\"obner basis for
 a family of multiobjective programs where the right-hand side varies.
 This new structure extends the notion of Gr\"obner basis for the single objective case, to the case of multiple objectives, i.e., a partial ordering
instead of a total ordering over the feasible vectors. The main
property of these bases is that the partial reduction
  of the integer elements in the kernel of the constraint matrix by the different blocks of the basis is zero.
  It allows us to prove that this new construction is a test family for a
   family of multiobjective programs. An algorithm '\`a la Buchberger' is
   developed to compute partial Gr\"obner bases and two different approaches are derived, using this methodology, for computing the entire set of efficient
solutions of any multiobjective integer linear problem (MOILP).
Some examples illustrate the application of the algorithms and
computational experiments are reported on several families of
problems.

 \end{abstract}

\section{Introduction}

The multiobjective paradigm appeared in economic theory in the
nineteenth century in the seminal works by Edgeworth
\cite{edgeworth1881} and Pareto \cite{pareto1896} to define an
economic equilibrium. Mathematically, the multiobjective
optimization approach consists of determining the maximal
(minimal) elements of a partially ordered set. This problem was
already addressed by Cantor \cite{cantor1897}, Cayley
\cite{cayley1889} and Hausdorff \cite{hausdorff1906} at the end of
the nineteenth century. Since then, multiobjective programming
(including multicriteria optimization) has been a fruitful
research field within the areas of applied mathematics, operations
research, and economic theory. Excellent textbooks and survey
papers are available in the literature, the interested reader is
referred to the books by Sawaragi, Nakayama and Tanino
\cite{tanino85}, Chankong and Haimes \cite{chankong-haimes83}, Yu
\cite{yu74}, Miettinen \cite{miettinen99} or Ehrgott, Figueira and
Gandibleux \cite{ehrgott06}, and to the surveys in
\cite{ehrgott02} and \cite{ehrgott05}.

The importance of multiobjective optimization is not only due to
its theoretical implications but also to its many applications.
Witnesses of that are the large number of real-world decision
problems that appear in the literature formulated as
multiobjective programs. Examples of them are flowshop scheduling
(see ~\cite{ishibuchi-murata98}), analysis in finance (see
\cite{ehrgott02}, Chapter 20), railway network infrastructure
capacity (see \cite{delorme03}), vehicle routing problems (see
\cite{jozefowiez04, sherbeny01}) or trajectory optimization (see
\cite{steuer85}) among many others.

Multiobjective programs are formulated as optimization (without
lost of generality, we restrict ourselves to the minimization
case) problems over feasible regions with at least two objective
functions. Usually, it is not possible to minimize all the
objective functions simultaneously since the objective functions
induce a partial order over the vectors in the feasible region, so
a different notion of solution is needed. A feasible vector is
said to be Pareto-optimal (efficient or non-dominated) if no other
feasible vector has componentwise smaller objective values, with
at least one strict inequality.

This paper studies multiobjective integer linear programs (MOILP).
Thus, we assume that all objective functions and constraints that
define the feasible region are linear, and that the feasible
vectors have non-negative integer components.

There are nowadays several exact methods to solve MOILP (see
\cite{ehrgott02}). Two of them claimed to be of general use and
have attracted the attention of researchers over the years:
multiobjective implicit enumeration (see \cite{zionts79},
\cite{zionts-wallenius80}) and multiobjective dynamic programming
(see \cite{karwan-villareal82}). Nevertheless, although in
principle they may be applied to any number of objectives, one can
mainly find, in the literature, applications to bicriteria
problems. On the other hand, there are several methods that apply
to bicriteria problems but that do not extend to the general case.
Thus, one can see that there are two thresholds in multiobjective
programming, a first step from 1 to 2 objectives and a second, and
deeper one, from 2 to more than two objectives. Thus, most of the
times, algorithms to solve multiobjective integer problems are
designed to compute only the solutions for the bicriteria case.
Moreover, some methods even do not provide the entire set of
Pareto-optimal solutions, but the supported ones (those that can
be obtained as solutions of linearly scalarized programs).

It is worth noting that most MOILP problems are NP-hard and
intractable (see \cite{ehrgott-gandibleux00} for further details).
Even in most cases where the single-objective problem is
polynomially solvable the multiobjective version becomes NP-hard.
This is the case of spanning tree problems and min-cost flow
problems, among others (see \cite{hamacher94} and
\cite{ehrgott00}). Therefore, computational efficiency is not an
issue when analyzing MOILP. The important point is to develop
tools that can handle these problems and that give insights into
their intrinsic nature. The goal of this paper is to present a new
general methodology for solving MOILP using tools borrowed from
algebraic geometry. The usage of algebraic geometry tools in
integer programming (single criterion) is not new (see
\cite{conti-traverso91}, \cite{hosten-sturmfels95},
\cite{thomas95}, \cite{hosten97}, \cite{urbaniak97},
\cite{thomas-weismantel97}). The main idea is to compute a
Gr\"obner basis for certain toric ideals (related to the
constraints matrix) with a monomial order induced by the objective
function.

Gr\"obner bases were introduced  by Bruno Buchberger in 1965 in
his PhD Thesis \cite{buchberger65}. He named them Gr\"obner bases
paying tribute to his advisor Wolfgang Gr\"obner. This theory
emerged as a generalization, from the one variable case to the
multivariate polynomial case, of the greatest common divisor in an
ideal sense. One of the outcomes of Gr\"obner bases theory was its
application to Integer Programming, firstly published by Conti and
Traverso ~\cite{conti-traverso91}. After this paper, a number of
publications using Gr\"obner bases to solve integers programs
appeared in the literature.

In ~\cite{hosten-sturmfels95}, Hosten and Sturmfels gave two ways
to implement Conti and Traverso algorithm that improve in many
cases branch-and-bound algorithm to solve, exactly, integer
programs. Thomas presented in ~\cite{thomas95} a geometric point
of view of the Buchberger algorithm as a method to obtain
solutions of an integer program. Later, Thomas and Weissmantel
\cite{thomas-weismantel97} improved the Buchberger algorithm in
its application to solve integer programs introducing truncated
Gr\"obner bases. At the same time, Urbaniak et al
\cite{urbaniak97} published a clear geometric interpretation of
the reduction steps of this kind of algorithms in the original
space (decision space). The interested reader can find excellent
descriptions of this methodology in the books by Adams and
Loustanau \cite{adams94}, Sturmfels \cite{sturmfels95}, Cox et al
\cite{cox-little-oshea98} or Bertsimas and Weissmantel
\cite{bertsimas05}, and in the papers by Aardal et al.
\cite{aardal02}, Sturmfels \cite{sturmfels02}, \cite{sturmfels04},
Sturmfels and Thomas \cite{sturmfels-thomas97} and Thomas
\cite{thomas97}

The main contribution of this paper is to adapt some tools from
algebraic geometry to solve multiobjective integer linear
programs. We present in this paper an algorithm to solve exactly
multiobjective problems, i.e. providing the whole set of
Pareto-optimal solutions (supported and non-supported ones). One
of the main advantages of our approach is that the number of
objective functions does not increase significantly the
difficulty. A new geometric approach of the concept of reduction
based on a partial ordering is given. This reduction allows us to
extend the concept of Gr\"obner basis when a partial ordering
rather than a total order is considered over $\N^n$. We call these
new structures partial Gr\"obner bases or p-Gr\"obner bases. We
prove that p-Gr\"obner bases can be generated by a variation of
the Buchberger algorithm in a finite number of steps. The main
property of a p-Gr\"obner basis is that, for each pair in
$\Z^n\times\Z^n_+$ with first component in $Ker(A)$, the reduction
by maximal chains in the basis is the zero set.

We propose two different algorithms to solve multiobjective
integer programs based on this new construction. Our first
algorithm consists of three stages. The first one only uses the
constraint matrix of the problem and it produces a system of
generators for the toric ideal $I_A$ (or its geometric
representation, $\Im_A$). In the second step, a p-Gr\"obner basis
is built using the initial basis given by the system of generators
computed in the first step. This step requires to fix the
objective matrix since it induces the partial order used in the
reduction steps. Once the right-hand side is fixed, in the third
step the Pareto-optimal solutions are obtained. This computation
uses the new concept of partial reduction of an initial feasible
solution by the p-Gr\"obner basis.

This algorithm extends, to some extent, Hosten-Sturmfels'
algorithm \cite{hosten-sturmfels95} for integer programs because
if we apply our method to single-objective problems, partial
reductions and p-Gr\"obner bases coincide with the standard
notions of reductions and Gr\"obner bases, respectively.

Our second algorithm is based on the original idea by Conti and
Traverso \cite{conti-traverso91}. It consists of using the big-M
method that results in an increasing number of variables, in order
to have an initial system of generators. Moreover, this approach
also provides an initial feasible solution. Therefore, the first
step in the above algorithm can be ignored and the third step is
highly simplified. In any case, our first algorithm (the one
extending Hosten-Sturmfels approach) has proved to be more
efficient than this second one since computation of a p-Gr\"obner
basis is highly sensitive to the number of variables.

Both algorithms have been implemented in MAPLE 10. In this paper
we report on some computational experiments based on two different
families of problems with different number of objective functions.

This paper is organized as follows. In Section \ref{sect1} we give
the notation, the formulation of the problem, and its algebraic
codification. In this section we also introduce the notion of test
family and its geometric description. Section \ref{sect2} presents
the definition of p-Gr\"obner basis, based on the notion of
partial reduction. Here, we also state the relationship between
test families and p-Gr\"obner bases: the reduced p-Gr\"obner basis
for a family of multiobjective programs varying the right-hand
side coincides with the minimal test family for that family. At
the end of the section, an illustrative example is presented.
Section \ref{sect3} is devoted to present the results of the
computational experiments and its analysis. Here, we solve several
families of MOILP, report on the performance of the algorithms and
draw some conclusions on their results and their implications.

\section{The problem and its translation}
\label{sect1}

The goal of this paper is to solve the multiobjective integer
linear program (MOILP) in its standard form:
\begin{eqnarray}
\label{MOILP_gral}
\min &\; (c^1\,x, \ldots, c^k\,x)\nonumber\\
s.t.&\nonumber\\
& \dsum_{j=1}^n \,a_{ij}\,x_j = b_i & i=1, \ldots, m\\
& x_j \in \Z_+ & j=1, \ldots, n\nonumber
\end{eqnarray}

with $b_i$ nonnegative integers, $x_i$ non negative and the
constraints are defining a polytope (bounded). Let us denote by
$A=(a_{ij}) \in \Z^{m\times n}$, $b=(b_{i}) \in \Z^m_+$ and
$C=(c_{ij}) \in \Z^{k\times n}_+$. In the following, Problem
\eqref{MOILP_gral} will be referred to as $MIP_{A,C}(b)$  and we
denote by $MIP_{A,C}$ the family of multiobjective problems where
the right-hand side varies.

The reader may note that there is no loss of generality in our
approach to multiobjective integer linear programming since any
general multiobjective integer linear problem with inequality
constraints and rational components in $A$, $b$ and $C$ can be
transformed to a problem in the above standard form.

It is clear that the problem $MIP_{A,C}(b)$ is not an usual
optimization problem since the objective function is a vector,
thus inducing a partial order among its feasible solutions. Hence,
solving the above problem requires an alternative concept of
solution, namely the set of non-dominated or Pareto-optimal points
(vectors).

A feasible vector $\widehat{x} \in \R^n$ is said to be a
\textit{Pareto-optimal solution} of $MIP_{A,C}(b)$ if there is no
other feasible vector $y$ such that
$$
c_j\,y \leq c_j\,\widehat{x} \qquad \forall j=1, \ldots, k
$$
with at least one strict inequality for some $j$.

If $x$ is a Pareto-optimal solution, the vector $(c_1\,x, \ldots,
c_k\,x)$ is called \textit{efficient}.

We say that a feasible point, $y$, is dominated by a feasible
point $x$ if $c_i\,x \leq c_i\,y$ for all $i=1, \ldots, k$, with
at least one strict inequality. According to the above concept,
solving a multiobjective problem consists of finding its entire
set of Pareto-optimal solutions, including those that have the
same objective values.

From the objective function $C$, we obtain a partial order over
$\Z^n$ as follows:

$$
x \prec_C y : \Longleftrightarrow C\,x \lvertneqq C\,y \quad or
\quad x=y
$$
where $Cx \lvertneqq Cy$ stands for $Cx \leq Cy$ and $Cx\neq Cy$.

Notice that since $C \in \Z_+^{m\times n}$, the above relation is
not complete. Hence, there may exist incomparable vectors (those
$x$, $y \in \Z^n_+$ such that neither $x\prec_C y$ or $y\prec_C
x$). We use this partial order, induced by the objective function
of Problem $MIP_{A,C}$ as the input for the multiobjective integer
programming algorithm developed in this paper.

\begin{remark}
The above order distinguishes solutions with the same objective
values and handles them as incomparable. This order can be refined
so that those solutions with the same objective values are not
incomparable. Consider the binary relation:
$$
x \preceq_C y : \Longleftrightarrow \left\{ \begin{array}{lr} C\,x
\lvertneqq C\,y & \text{ or}\\
Cx = Cy \text{ and } x \prec_{lex} y&
\end{array}\right.
$$
This alternative order allows us to rank those solutions that have
the same objective values using the lexicographical order of their
components.

Let us consider the following equivalence relation in $\Z^n$:
$$
x \sim_C y :\Longleftrightarrow Cx = Cy
$$
The above partial order, $\preceq_C$, allows us to solve a
simplified version of the multiobjective problem. In this version,
we obtain solutions in $\Z^n/\sim_C$, where $x \sim_C y
:\Longleftrightarrow Cx = Cy$. The reader may note that when
solving the problem with the order $\preceq_C$, one would obtain
only a representative element of each class of Pareto-optimal
solutions (the lexicographically smallest). With those efficient
values, $\{v_1, \ldots, v_t\}$, the remaining solutions can be
obtained solving the following system of diophantine equations, in
$x$, for each $v_i$, $i=1, \ldots, t$:
$$
\left\{\begin{array}{ll}
Cx &= v_i\\
Ax &= b\\
x &\in \Z^n_+
\end{array}\right.
$$
\end{remark}

\begin{remark}
\label{remark:slacks} In some cases, the order $\prec_C$ can be
refined to be adapted to specific problems. This is the case when
slack variables appear in mathematical programs. Two feasible
solutions $(x, s_1)$ and $(x, s_2)$, where $s_1$ and $s_2$ are the
slack components, have the same objective values. The order
$\prec_C$ considers both solutions as incomparable, although they
are the same because we are looking
 just for the $x$-part of the solution. In these cases, we consider the following refined partial order in $\Z^n \times
 \Z^r$,
$$
(x,s) \prec^s_C (y,s^{'}) : \Longleftrightarrow \left\{
\begin{array}{lr} C\,x
\lvertneqq C\,y & \text{ or}\\
Cx = Cy \text{ and } s \prec_{lex} s^{'}&
\end{array}\right.
$$
where $x, y \in \Z^n_+$ are the actual decision variables and $s,
s^{'} \in \Z^r_+$ the slack variables of our problem.
\end{remark}
In the following we will use the partial order $\prec_C$ unless it
is explicitly specified.

Our matrix $A$ is encoded in the set
$$
I_A = \{ \{u,v\} : u, v \in \N^n , u - v \in Ker(A)\}.
$$
Let $\pi: \N^n \longrightarrow \Z^n$ denote the map $x \mapsto
Ax$. Given a right-hand side vector $b$ in $\Z^n$, the set of
feasible solutions to $MIP_{A,C}(b)$ constitutes $\pi^{-1}(b)$,
the preimage of $b$ under this map. In the rest of this paper, we
identify the discrete set of points $\pi^{-1}(b)$ with its convex
hull and we call it the $b$-fiber of $MIP_{A,C}$. Thus,
$\pi^{-1}(b)$ or the $b$-fiber of $MIP_{A,C}$ is the polyhedron
that is the convex hull of all feasible solutions to
$MIP_{A,C}(b)$.

For any pair $\{u, v\}$, with $u, v\in \N^n$, we define the set
$setlm(u,v)$ as follows:
$$
setlm(u,v) = \left\{ \begin{array}{ll} \{u\} & \mbox{if $v \prec_C u$}\\
\{v\} & \mbox{if $u \prec_C v$}\\
\{u,v\} & \mbox{if $u$ and $v$ are incomparable by $\prec_C$}
\end{array}\right.
$$

The reader may note that $setlm(u,v)$ is the set of degrees of the
leading monomials according to the identification $\{u,v\} \mapsto
x^u-x^v \in \R[x_1, \ldots, x_n]$, induced by the partial order
$\prec_C$.

From the above definition, $setlm(u,v)$ may have more than one
leading term, since $\prec_C$ is only a partial order. To account
for all this information we denote by $\mathcal{F}(u,v)$ the set
of triplets
$$
\mathcal{F}(u,v) = \{ (u,v,w) : w \in setlm(u,v)\}.
$$
The above concept extends to any finite set of pairs of vectors in
$\N^n$, accordingly. For a pair of sets $\mathbf{u} = \{u_1,
\ldots, u_t\}$ and $\mathbf{v} =\{v_1,\ldots, v_t\}$ the
corresponding set of ordered pairs is:
$$
\mathcal{F}(\mathbf{u},\mathbf{v}) = \{ (u_i,v_i,w) : w \in
setlm(u_i,v_i), i=1, \ldots,t \}.
$$

$\mathcal{F}(\mathbf{u},\mathbf{v})$ can be partially ordered
based on the third component of its elements. Therefore, we can
see $\mathcal{F}(\mathbf{u},\mathbf{v})$ as a directed graph
$G(E,V)$ where $V$ is identified with the elements of
$\mathcal{F}(\mathbf{u},\mathbf{v})$ and
$((u_i,v_i,w),(u_j,v_j,w^{'})) \in E$ if
$(u_i,v_i,w),(u_j,v_j,w^{'}) \in V$ and $w^{'} \prec_C w$. We are
interested in the maximal ordered chains of $G$. Note that they
can be efficiently computed by different methods, see e.g.
~\cite{baer69}, ~\cite{schrijver03}.

 The above concepts are clarified in the next example.

\begin{ej}
\label{ej1}
Let $\mathbf{u}$ $=$ $\{(2,3)$,$(0,2)$,$(3,0)$,$(2,1)$,$(1,1)\}$,
 $\mathbf{v}$ $=$ $\{(1,4)$,$(1,3)$,$(4,2)$,$(1,2)$,$(1,0)\}$ and $\prec_C$ the partial order induced by the matrix
$$
C = \left[ \begin{array}{cc} 2 & 1 \\ 3 & 5 \end{array} \right]
$$
then, $setlm((2,3),(1,4))$ $=$ $\{(2,3)$,$(1,4)\}$,
$setlm((0,2),(1,3))$ $=$ $\{(1,3)\}$, $setlm((3,0),(4,2))$ $=$
$\{(4,2)\}$,\\
$setlm((2,1),(1,2))$ $=$ $\{(2,1)$,$(1,2)\}$ and
$setlm((1,1),(1,0))$ $=$ $\{(1,1)\}$. Now, by definition we have:
\begin{eqnarray*}
\begin{split}
\mathcal{F}(\mathbf{u},\mathbf{v}) = &\{ &\big((2,3), (1,4), (2,3)\big), &\big((2,3), (1,4), (1,4)\big), &\big((0,2), (1,3), (1,3)\big), \\
& &\big((3,0), (4,2), (4,2)\big), &\big((2,1), (1,2), (2,1)), &((2,1), (1,2), (1,2)\big), \\
& &\big((1,1), (1,0), (1,1)\big)  &\}. &
\end{split}
\end{eqnarray*}
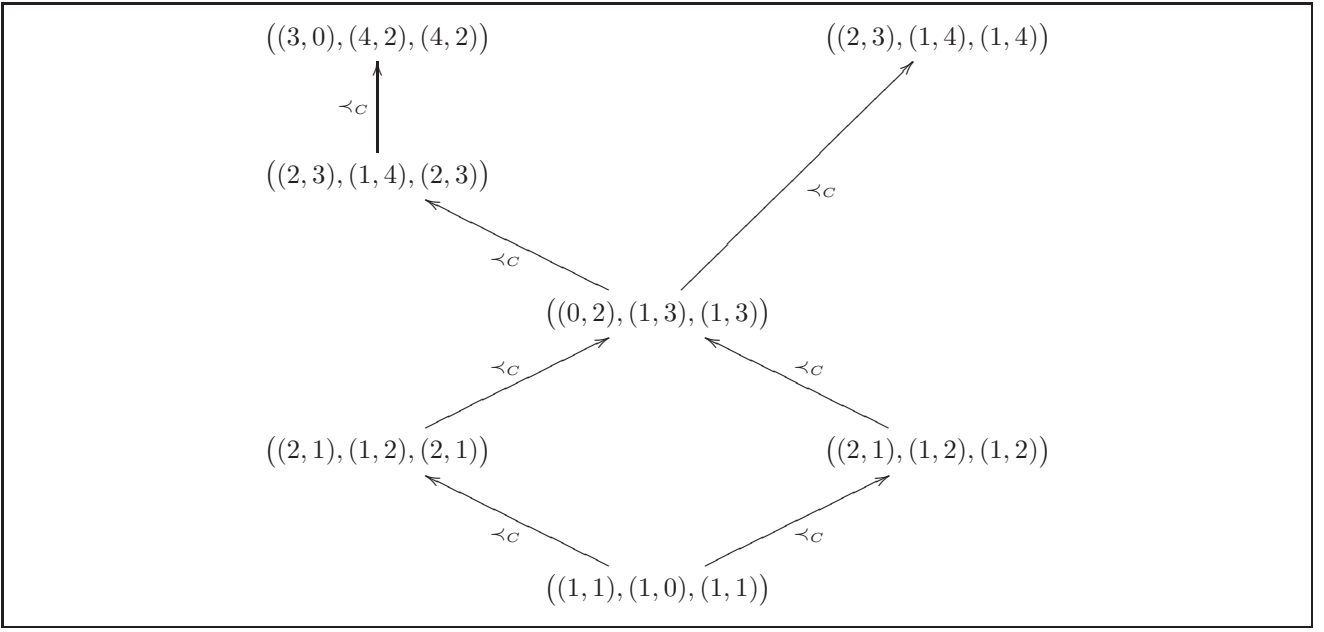
\begin{figure}

\begin{center}
\fbox{\parbox[t][\height][t]{\textwidth}{
$$
\xymatrix@C-10pt@R+10pt{\txt{$\big((3,0),(4,2),(4,2)\big)$} &  &\txt{$\big((2,3),(1,4), (1,4)\big)$}\\
\txt{$\big((2,3),(1,4), (2,3)\big)$}\ar[u]^{\prec_C}& & \\
&\txt{$\big((0,2),(1,3), (1,3)\big)$}\ar[uur]_{\prec_C} \ar[ul]^{\prec_C}&\\
\txt{$\big((2,1),(1,2),(2,1)\big)$}\ar[ur]^{\prec_C}& & \txt{$\big((2,1),(1,2),(1,2)\big)$}\ar[ul]_{\prec_C}\\
& \txt{$\big((1,1),(1,0),(1,1)\big)$} \ar[ur]_{\prec_C}
\ar[ul]^{\prec_C}&}
$$
}}
\end{center}
\caption{\label{hasse}Hasse diagram of the graph associated with
the data in Example \ref{ej1}}
\end{figure}

Figure \ref{hasse} corresponds to the directed graph associated
with $\mathcal{F}(\mathbf{u},\mathbf{v})$, according to the
partial ordering induced by $C$. There are four maximal chains:

$M_1= \{ \big((3,0),(4,2),(4,2)\big), \big((2,3),(1,4), (2,3)\big), \big((0,2),(1,3), (1,3)\big), \big((2,1),(1,2),(2,1)\big), \big((1,1),(1,0),(1,1)\big)\}$\\
$M_2= \{ \big((3,0),(4,2),(4,2)\big), \big((2,3),(1,4), (2,3)\big), \big((0,2),(1,3), (1,3)\big), \big((2,1),(1,2),(1,2)\big), \big((1,1),(1,0),(1,1)\big)\}$\\
$M_3= \{\big((2,3), (1,4), (1,4)\big), \big((0,2), (1,3), (1,3)\big), \big((2,1), (1,2), (2,1)\big), \big((1,1), (1,0), (1,1)\big)\}$\\
$M_4= \{\big((2,3),(1,4), (1,4)\big), \big((0,2),(1,3),
(1,3)\big), \big((2,1),(1,2),(1,2)), ((1,1),(1,0),(1,1)\big)\}$.
\end{ej}
For any pair of sets $\mathbf{u} = \{u_1, \ldots, u_t\}$ and
$\mathbf{v} =\{v_1,\ldots, v_t\}$ with $\{u_i, v_i\} \in I_A$, the
corresponding set $\mathcal{F}(\mathbf{u},\mathbf{v})$ may be seen
as a set of pairs in $\Z^n \times \Z^n_+$ through the following
map

\begin{eqnarray*}
\begin{aligned}
\phi : & \N^n \times \N^n \times \N^n &\longrightarrow & \Z^n \times \Z^n_+\\
& (u, v, w) & \mapsto & (u - v , w).
\end{aligned}
\end{eqnarray*}

Then, the maximal chains, $F_1, \ldots, F_t$, of the image of
$\mathcal{F}(\mathbf{u},\mathbf{v})$ under $\phi$ with respect to
the order $\prec_C$ over the second components, clearly satisfy
the following properties:
\begin{enumerate}
\item $F_i$ is totally ordered by the second components of its images via $\phi$ with respect to $\prec_C$, for $i=1, \ldots, t$.
\item For all $(\alpha, \beta) \in F_i$, $i=1, \ldots, t$, $A\,(\beta-\alpha) = A\,\beta$.
\end{enumerate}
The application $\phi$ and the above properties allow us to define
the notion of test family for $MIP_{A,C}$. This notion is
analogous to the concept of test set for a family of single
objective integer programs when we have a partial order rather
than a total order over $\N^n$ (see  \cite{thomas95}). Test
families are instrumental for finding the Pareto-optimal set of
each member $MIP_{A,C}(b)$ of the family of multiobjective integer
linear programs.

\begin{defi}[Test Family]\label{def:testfam}
A finite collection $\grob = \{\grob_C^1, \ldots, \grob_C^r\}$ of
sets in $\Z^n \times \Z^n_+$ is a test family for $MIP_{A,C}$ if
and only if:
\begin{enumerate}
\item $\grob_C^j$ is totally ordered by the second component with respect to $\prec_C$, for $j=1, \ldots, r$.
\item For all $(g,h) \in \grob_C^j$, $j=1, \ldots, r$, $A\,(h-g) = A\,h$.
\item If $x \in \N^n$ is a dominated solution for $MIP_{A,C}(b)$, with $b \in \Z^n_+$, there is some $\grob_C^j$ in the collection and $(g,h) \in  \grob_C^j$, such that $x - g \prec_C x$.
\item If $x \in \N^n$ is a Pareto-optimal solution for $MIP_{A,C}(b)$, with $b \in \Z^n_+$, then for all $(g,h) \in \grob_C^j$ and for all $j=1, \ldots, n$ either $x-g$ is infeasible or $x-g$ does not compare with $x$.
\end{enumerate}
\end{defi}

Given a test family for $MIP_{A,C}$ there is a natural approach to
find the entire Pareto-optimal set. Suppose we wish to solve
$MIP_{A,C}(b)$ for which $x^*$ is a feasible solution.

If $x^*$ is dominated then there is some $j$ and $(g,h) \in
\grob_C^j$ such that $x^* - g$ is feasible and $x^* -g \prec_C
x^*$, whereas for the remaining chains there may exist some
$(g,h)$ such that $x^*-g$ is feasible but incomparable with $x^*$.
We keep tracks of all of them.

If $x^*$ is non-dominated, we have to keep it as an element in our
current solution set. Then, reducing $x^*$ by the chains in the
test family we can only obtain either incomparable feasible
solutions, that we maintain in our structure, or infeasible
solutions that are discarded.

The above two cases lead us to generate the following set. From
$x^*$ we compute the set of incumbent solutions:
\begin{center}
$ IS(x^*) := \{y^* : y^* = x^*-g_{j_i}, (g_{j_i}, h_{j_i}) \mbox{ is the largest element $(g,h)$ in the chain}$ \\
$\grob_{C}^i$ such that $x^* - g$ is feasible , $i=1, \ldots,
r\}$.
\end{center}

Now, the scheme proceeds recursively on each element of the set
$IS(x^*)$. Finiteness of the above scheme is clear since we are
generating a search tree with bounded depth (cardinality of the
test family) and bounded width, each element in the tree has at
most $r$ (number of chains) followers. Correctness of this
approach is ensured since any pair of non-dominated solutions must
be connected by a reduction chain through elements in the test
family (see Theorem \ref{theo1} and Corollary \ref{coro1}).

The above approach assumes that a feasible solution to
$MIP_{A,C}(b)$ is known (thus implying that the problem is
feasible). Methods to detect infeasibility and to get an initial
feasible solution are connected to solving diophantine systems of
linear equations, the interested reader is referred to
\cite{pottier91}, for further details.

The following lemmas help us in describing the geometric structure
of a test family for multiobjective integer linear problems.

\begin{lemma}[Gordan-Dickson Lemma, Theorem 5 in \cite{cox-little-oshea92}]
\label{lemma:gordan-dickson} If $P \subseteq \N^n$, $P \neq
\emptyset$, then there exists a minimal subset $\{p_1, \ldots,
p_m\} \subseteq P$ that is finite and unique such that $p \in P$
implies $p_j \leq p$ (component-wise) for at least one $j=1,
\ldots, m$.
\end{lemma}

\begin{lemma}\label{lemma_thomas} There
exists a unique, minimal, finite set of vectors $\alpha_1, \ldots,
\alpha_k \in \N^n$ such that the set $\mathcal{L}_C$ of all
dominated solutions in all fibers of $MIP_{A,C}$ is a subset of
$\N^n$ of the form
$$
\mathcal{L}_{C} = \bigcup_{j=1}^{k} (\alpha_j + \N^n).
$$
\end{lemma}
\begin{proof} The set of dominated solutions of all
problems $MIP_{A,C}$ is:
$$
\mathcal{L}_C = \{ \alpha \in \N^n : \exists \beta \in \N^n \mbox{
with } A\beta = A\alpha \mbox{ and } \beta \prec_C \alpha\}.
$$
Let $\alpha$ be an element in $\mathcal{L}_C$ and $\beta$ a
Pareto-optimal point in the fiber $\pi^{-1}(A\alpha)$ that
satisfies $\beta \prec_C \alpha$. Then, for any $\gamma \in \N^n$,
$A(\alpha+\gamma)=A(\beta+\gamma)$, $\alpha+\gamma, \beta+\gamma
\in \N^n$ and $\beta+\gamma \prec_C \alpha+\gamma$, because the
cost matrix, $C$, has only nonnegative coefficients. Therefore,
$\alpha + \gamma$ is a feasible solution dominated by
$\beta+\gamma$ in the fiber $\pi^{-1}(A(\alpha+\gamma))$. Then,
$\alpha + \gamma \in \mathcal{L}_C$ for all $\gamma \in \N^n$, so,
$\alpha + \N^n \subseteq \mathcal{L}_C$. By Lemma
\ref{lemma:gordan-dickson} we conclude that there exists a minimal
set of elements $\alpha_1, \ldots, \alpha_k \in \N^n$ such that
$\mathcal{L}_C = \bigcup_{j=1}^{k} (\alpha_j + \N^n)$.
\end{proof}

 Once the elements $\alpha_1,
\ldots, \alpha_k$ that generates $\mathcal{L}_C$ (in the sense of
the above result) have been obtained, one can compute the maximal
chains of the set $\{\alpha_1, \ldots, \alpha_k\}$ with respect to
the partial order $\prec_C$. We denote by $\mathcal{C}_C^1,
\ldots, \mathcal{C}_C^\mu$ these maximal chains and set
$\mathcal{L}_C^i = \bigcup_{t=1}^{k_i} (\alpha_t^i + \N^n)$, where
$\alpha_t^i \in \mathcal{C}_C^i$ for $t=1, \ldots,k_i$ and $i=1,
\ldots, \mu$. For details about maximal chains, upper bounds for
its cardinality and algorithms to compute them for a partially
ordered set, the reader is refereed to \cite{baer69}.

It is clear that, with this construction, we have: $\mathcal{L}_C
= \dbigcup_{i=1}^\mu \mathcal{L}^i_C$.

We now describe a finite family of sets $\grob_{\prec_C} \subseteq
Ker(A) \cap \Z^n$ and prove that it is indeed a test family for
$MIP_{A,C}$.

Let $\grob_{\prec_C} = \{ \grob^i_{\prec_C} \}_{i=1}^{\mu}$, where
\begin{equation}
\label{gc} \grob_{\prec_C}^i= \{({g}^k_{ij},{h}^k_{ij}) =
(\alpha^i_j-\beta_{ij}^k, \alpha^i_j) , j=1, \ldots k_i, k=1,
\ldots, m_{ij}\}, i=1, \ldots,\mu,
\end{equation}
are the maximal chains of $\grob_{\prec_C}$ (with respect to the
order $\prec_C$ over the second components) and where $\alpha^i_1,
\ldots, \alpha^i_{k_i}$ are the unique minimal elements of
$\mathcal{L}^i_{\prec_C}$ and $\beta_{ij}^1, \ldots,
\beta_{ij}^{m_{ij}}$ the Pareto-optimal solutions to the problem
$MIP_{A,C}(A\alpha^i_j)$.

In the next section we give an algorithm that explicitly
constructs $\grob_{\prec_C}$. Notice that for fixed $i, j$ and
$k$, $g_{ij}^k=(\alpha^i_j - \beta_{ij}^k)$ is a point in the
subspace $S=\{ x \in \Q^n : Ax =  0\}$, i.e., in the 0-fiber of
$MIP_{A,C}$. Geometrically we think of $(\alpha^i_j -
\beta_{ij}^k, \alpha^i_j)$ as the oriented vector
$\overrightarrow{g}^k_{ij} = \overrightarrow{[\alpha^i_j,
\beta^k_{ij}]}$ in the $A\alpha^i_j$-fiber of $MIP_{A,C}$ directed
to the Pareto-optimal solution $\beta^k_{ij}$. The vector is
directed from the non-optimal point $\alpha^i_{j}$, to the
Pareto-optimal point $\beta^k_{ij}$ due to the minimization
criterion in $MIP_{A,C}$ which requires us to move away from
expensive points. Subtracting the point $\overrightarrow{g}^k_{ij}
= \alpha^i_j - \beta^k_{ij}$ to the feasible solution $\gamma$
gives the new solution $\gamma - \alpha^i_j + \beta^k_{ij}$ which
is equivalent to translating $\overrightarrow{g}^k_{ij}$ by a
nonnegative integer vector.

Consider an arbitrary fiber of $MIP_{A,C}$ and a feasible lattice
point $\gamma$ in this fiber. For each vector
$\overrightarrow{g}^k_{ij}$ in $\grob_{\prec_C}$, check whether
$\gamma - g_{ij}^k$ is in $\N^n$. At $\gamma$ draw all such
possible translations of vectors from $\grob_{\prec_C}$. The head
of the translated vector is also incident at a feasible point in
the same fiber as $\gamma$ since $g_{ij}^k$ is in the 0-fiber of
$MIP_{A,C}$. We do this construction for all feasible points in
all fibers of $MIP_{A,C}$. From Lemma \ref{lemma_thomas} and the
definition of $\grob_{\prec_C}$, it follows that no vector in
$\grob_{\prec_C}$ can be translated by a $\nu$ in $\N^n$ such that
its tail meets a Pareto-optimal solution on a fiber unless the
obtained vector is incomparable with the Pareto-optimal point.

\begin{theo}\label{theo1}
The above construction builds a connected directed graph in every
fiber of $MIP_{A,C}$. The nodes of the graph are all the lattice
points in the fiber and $(\gamma, \gamma^{'})$ is an edge of the
directed graph if $\gamma^{'} = \gamma - g_{ij}^k$ for some $i$,
$j$ and $k$. For each maximal chain in the $b$-fiber of
$MIP_{A,C}$, its directed graph has a unique final node at each
Pareto-optimal solution for $MIP_{A,C}(b)$.
\end{theo}

\begin{proof}
Pick a fiber of $MIP_{A,C}$ and at each feasible lattice point
construct all possible translations of the vector
$\overrightarrow{g}^k_{ij}$ from the set $\grob_{\prec_C}^i$ as
described above. Let $\alpha$ be a lattice point in this fiber. By
Lemma \ref{lemma_thomas}, $\alpha = \alpha^i_j + \nu$ for some $i
\in \{1, \ldots, t\}$ and $\nu \in \Z_+^n$. Now, since
$\alpha_k^{'} = \beta^k_{ij} + \nu$ also lies in this fiber, then
$\alpha^{'}_k \prec_C \alpha$ or $\alpha^{'}_k$ and $\alpha$ are
incomparable. Therefore, $\overrightarrow{g}^k_{ij}$ translated by
$\nu \in \N^n$ is an edge of this graph and we can move along it
from $\alpha$ to a point $\alpha^{'}$ in the same fiber, such that
$\alpha^{'}\prec_C \alpha$ or $\alpha$ and $\alpha^{'}$ are
incomparable. This proves that from every dominated point in the
fiber we can reach an improved or incomparable point (with respect
to $\prec_C$) in the same fiber by moving along an edge of the
graph.

By the construction above, the outdegree of any terminal element
in any maximal chain is $0$. Therefore, any directed maximal path
from a dominated point must end exactly at one Pareto-optimal
point.
\end{proof}

We call the graph in the $b$-fiber of $MIP_{A,C}$ built from elements in $\grob_{\prec_C}$, the $\prec_C$-skeleton of that fiber.

The reader may note that from each dominated solution $\alpha$,
one can easily build paths to its comparable Pareto-optimal
solutions subtracting elements in $\grob_{\prec_C}$. Indeed, let
$\alpha_i$ be a minimal element of $\mathcal{L}_C$ such that
$\alpha = \alpha_i + \gamma$, with $\gamma \in \N^n$, and let
$\beta_i$ be the Pareto-optimal solution in the $A\alpha_i$-fiber
that is comparable with $\alpha_i$ and such that $\beta_i+\gamma$
is comparable with $\beta$. Then $\alpha^{'}=\beta_i+\gamma$ is a
solution in the $A\alpha$-fiber with $\beta \prec_C
\alpha^{'}\prec_C \alpha$. Now, one repeats this process but
starting with $\alpha^{'}$ and $\beta$, until $\alpha^{'}=\beta$.
Moreover, the case where $\alpha$ and $\beta$ are incomparable
reduces to the previous one by finding a path from $\alpha$ to any
intermediate point $\beta^{'}$ that compares with $\beta$. This
analysis leads us to the following result.

\begin{corollary}
\label{coro1} In the  $\prec_C$-skeleton of a fiber there exists a
directed path from every feasible point $\alpha$ to each
Pareto-optimal point, $\beta$, in the same fiber. The vectors of
objective function values of successive points in the path do not
increase componentwise from $\alpha$ to $\beta$.
\end{corollary}

\begin{corollary}\label{grobC} The
family $\grob_{\prec_C}$ is the unique minimal test family for
$MIP_{A,C}$. It depends only on the matrix $A$ and the cost matrix
$C$.
\end{corollary}

\begin{proof}
By definition of $\grob_{\prec_C}$, the conditions 1. and 2. of
Definition \ref{def:testfam} are satisfied. From Theorem
\ref{theo1} it follows that properties 3. and 4. are also
satisfied, so $\grob_{\prec_C}$ is a test family for $MIP_{A,C}$.
Minimality is due to the fact that removing any element
$({g}^k_{ij},{h}^k_{ij})$ from $\grob_{\prec_C}$ results in
$\grob_{\prec_C} \setminus \{ ({g}^k_{ij},{h}^k_{ij})\}$. However,
this new set is not a test family since no oriented vector in
$\grob_{\prec_C} \setminus \{ ({g}^k_{ij},{h}^k_{ij})\}$ can be
translated through a nonnegative vector in $\N^n$ such that its
tail meets $\alpha^i_j$. It is clear by definition that
$\grob_{\prec_C}$ depends only on $A$ and $C$.
\end{proof}

\begin{ej}
Let $MIP_{A,C}$ be the family of multiobjective problems, with the
following constraints and objective function matrices: \label{ej2}
$$
A=\left[ \begin{array}{cccc} 2 & 2 & -1 & 0\\ 0 & 2 & 0 & 1
\end{array}\right] ,\qquad C=\left[\begin{array}{cccc} 10 & 1 & 0 &
0\\ 1 & 10 & 0 & 0 \end{array}\right].
$$
Let $(x_1, x_2, s_1, s_2)$ be the vector of variables, where $s_1$
and $s_2$ are slack variables. In this example, using the order
$\prec^s_C$ (see Remark \ref{remark:slacks}), $\grob_C =
\{\grob_C^1, \grob_C^2\}$, where
$\grob_{\prec_C}^1=\{\overrightarrow{g}_1^1=((0,1,2,-1),(0,1,2,0)),\\
\overrightarrow{g}_2^1=((-1,1,0,-2),(0,1,0,0))\}$ and
$\grob_{\prec_C}^2=\{\overrightarrow{g}_1^2=((1,0,2,0),(1,0,2,0)),
\overrightarrow{g}_2^2=((1,-1,0,2),(1,0,0,2))\}$.

Figure \ref{fig1} shows, on the $(x_1, x_2)$-plane, the
$\prec_C$-skeleton of the fiber corresponding to the right-hand
side vector $(17,11)^t$. In the box over the graph of the
$\prec_C$-skeleton, we show the second components of the elements
of $\grob_{\prec_C}$. The reader may note that in the graph, the
arrows have opposite directions due to the fact that the directed
paths (improving solutions) are built subtracting the elements in
$\grob_{\prec_C}$. We describe how to compute the sets
$\grob_{\prec_C}^1$ and $\grob_{\prec_C}^2$ in Section
\ref{sect2}.

\begin{figure}

\center
\fbox{
\parbox[h]{0.6\textwidth}{
$
\xymatrix{
&& &  &\ar[rrdd]^{g_1^2} & & & & &\\
&&& & & & &  \ar[rr]_{g_2^2} & &\\
&& \ar@{.>}[uull]^{g_1^1} & \ar@{.>}[uu]^{g_2^1} & & & &&&}
$
\center\small
Elements in $\mathcal{G}_{\prec_C} = \{\mathcal{G}_{\prec_C}^1, \mathcal{G}_{\prec_C}^2\}$
}}
$$
\xy
\POS (80,50) *\cir<1pt>{}="a1"
, \POS (70,50) *\cir<1pt>{}="a2"
, \POS (60,50) *\cir<1pt>{}="a3"
, \POS (50,50) *\cir<1pt>{}="a4"
, \POS (40,50) *\cir<2pt>{}="a5"
, \POS (90,40) *\cir<1pt>{}="a6"
, \POS (80,40) *\cir<1pt>{}="a7"
, \POS (70,40) *\cir<1pt>{}="a8"
, \POS (60,40) *\cir<1pt>{}="a9"
, \POS (50,40) *\cir<2pt>{}="a10"
, \POS (100,30) *\cir<1pt>{}="a11"
, \POS (90,30) *\cir<1pt>{}="a12"
, \POS (80,30) *\cir<1pt>{}="a13"
, \POS (70,30) *\cir<1pt>{}="a14"
, \POS (60,30) *\cir<2pt>{}="a15"
, \POS (110,20) *\cir<1pt>{}="a16"
, \POS (100,20) *\cir<1pt>{}="a17"
, \POS (90,20) *\cir<1pt>{}="a18"
, \POS (80,20) *\cir<1pt>{}="a19"
, \POS (70,20) *\cir<2pt>{}="a20"
, \POS (120,10) *\cir<1pt>{}="a21"
, \POS (110,10) *\cir<1pt>{}="a22"
, \POS (100,10) *\cir<1pt>{}="a23"
, \POS (90,10) *\cir<1pt>{}="a24"
, \POS (80,10) *\cir<2pt>{}="a25"
, \POS (130,0) *\cir<1pt>{}="a26"
, \POS (120,0) *\cir<1pt>{}="a27"
, \POS (110,0) *\cir<1pt>{}="a28"
, \POS (100,0) *\cir<1pt>{}="a29"
, \POS (90,0) *\cir<2pt>{}="a30"

\POS"a1" \ar "a2"
\POS"a2" \ar "a3"
\POS"a3" \ar "a4"
\POS"a4" \ar^(1.25){{\text{\tiny{$(4,5,1,1)$}}}} "a5"
\POS"a6" \ar "a7"
\POS"a7" \ar "a8"
\POS"a8" \ar "a9"
\POS"a9" \ar^(1.25){{\text{\tiny{$(5,4,1,3)$}}}} "a10"
\POS"a11" \ar "a12"
\POS"a12" \ar "a13"
\POS"a13" \ar "a14"
\POS"a14" \ar^(1.25){{\text{\tiny{$(6,3,1,5)$}}}} "a15"
\POS"a16" \ar "a17"
\POS"a17" \ar"a18"
\POS"a18" \ar "a19"
\POS"a19" \ar^(1.25){{\text{\tiny{$(7,2,1,7)$}}}} "a20"
\POS"a21" \ar "a22"
\POS"a22" \ar "a23"
\POS"a23" \ar "a24"
\POS"a24" \ar^(1.25){{\text{\tiny{$(8,1,1,9)$}}}}"a25"
\POS"a26" \ar "a27"
\POS"a27" \ar "a28"
\POS"a28" \ar "a29"
\POS"a29" \ar^(1.25){{\text{\tiny{$(9,0,1,11)$}}}} "a30"
\POS"a26" \ar@<0.3ex> "a21"
\POS"a21" \ar@<0.3ex> "a16"
\POS"a16" \ar@<0.3ex> "a11"
\POS"a11" \ar@<0.3ex> "a6"
\POS"a6" \ar@<0.3ex> "a1"
\POS"a21" \ar@<0.3ex>@{.>} "a26"
\POS"a16" \ar@<0.3ex>@{.>} "a21"
\POS"a11" \ar@<0.3ex>@{.>} "a16"
\POS"a6" \ar@<0.3ex>@{.>} "a11"
\POS"a7" \ar@<0.3ex> "a2"
\POS"a12" \ar@<0.3ex> "a7"
\POS"a17" \ar@<0.3ex> "a12"
\POS"a22" \ar@<0.3ex> "a17"
\POS"a27" \ar@<0.3ex> "a22"
\POS"a12" \ar@<0.3ex>@{.>} "a17"
\POS"a17" \ar@<0.3ex>@{.>} "a22"
\POS"a22" \ar@<0.3ex>@{.>} "a27"
\POS"a8" \ar@<0.3ex> "a3"
\POS"a13" \ar\ar@<0.3ex> "a8"
\POS"a18" \ar\ar@<0.3ex> "a13"
\POS"a23" \ar@<0.3ex> "a18"
\POS"a28" \ar@<0.3ex> "a23"
\POS"a18" \ar@<0.3ex>@{.>} "a23"
\POS"a23" \ar@<0.3ex>@{.>} "a28"
\POS"a9" \ar@<0.3ex> "a4"
\POS"a14" \ar@<0.3ex> "a9"
\POS"a19" \ar@<0.3ex> "a14"
\POS"a24" \ar@<0.3ex> "a19"
\POS"a29" \ar@<0.3ex> "a24"
\POS"a24" \ar@<0.3ex>@{.>} "a29"
\POS"a10" \ar@<0.3ex> "a5"
\POS"a15" \ar@<0.3ex> "a10"
\POS"a20" \ar@<0.3ex> "a15"
\POS"a25" \ar@<0.3ex> "a20"
\POS"a30" \ar@<0.3ex> "a25"
\POS"a6" \ar@{.>} "a12"
\POS"a12" \ar@{.>} "a18"
\POS"a18" \ar@{.>} "a24"
\POS"a24" \ar@{.>} "a30"
\POS"a11" \ar@{.>} "a17"
\POS"a17" \ar@{.>} "a23"
\POS"a23" \ar@{.>} "a29"
\POS"a16" \ar@{.>} "a22"
\POS"a22" \ar@{.>} "a28"
\POS"a21" \ar@{.>} "a27"

\POS"a5" \ar@<0.3ex>@{.>} "a10"
\POS"a10" \ar@<0.3ex>@{.>} "a15"
\POS"a15" \ar@<0.3ex>@{.>} "a20"
\POS"a20" \ar@<0.3ex>@{.>} "a25"
\POS"a25" \ar@<0.3ex>@{.>} "a30"
\POS"a4" \ar@<0.3ex>@{.>} "a9"
\POS"a9" \ar@<0.3ex>@{.>} "a14"
\POS"a14" \ar@<0.3ex>@{.>} "a19"
\POS"a19" \ar@<0.3ex>@{.>} "a24"
\POS"a3" \ar@<0.3ex>@{.>} "a8"
\POS"a8" \ar@<0.3ex>@{.>} "a13"
\POS"a13" \ar@<0.3ex>@{.>} "a18"
\POS"a2" \ar@<0.3ex>@{.>} "a7"
\POS"a7" \ar@<0.3ex>@{.>} "a12"
\POS"a1" \ar@<0.3ex>@{.>} "a6"
\POS"a22" \ar@<0.3ex>@{.>} "a27"

\POS"a4" \ar@{.>} "a10"
\POS"a3" \ar@{.>} "a9"
\POS"a9" \ar@{.>} "a15"
\POS"a2" \ar@{.>} "a8"
\POS"a8" \ar@{.>} "a14"
\POS"a14" \ar@{.>} "a20"
\POS"a1" \ar@{.>} "a7"
\POS"a7" \ar@{.>} "a13"
\POS"a13" \ar@{.>} "a19"
\POS"a19" \ar@{.>} "a25"
\endxy
$$

\caption{\label{fig1}The $\prec_C$-skeleton  of the
$(17,11)^t$-fiber of $MIP_{A,C}$ projected on the $x_1,
x_2$-plane.}
\end{figure}
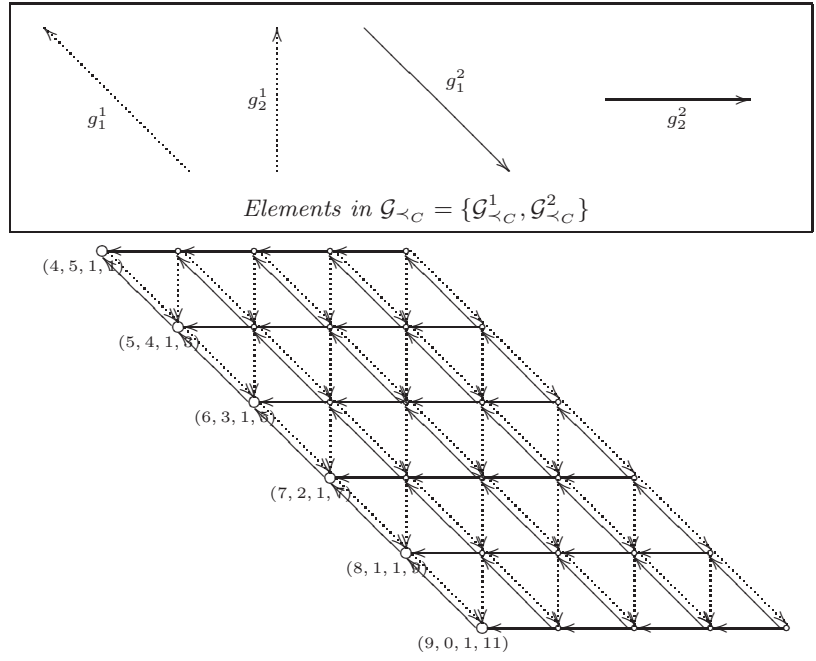

\end{ej}

Given $\grob_{\prec_C}$, there are several ways to build a path
from each feasible point in a fixed fiber to any Pareto-optimal
solution. However, there is a canonical way to do it: Fix $\sigma$
a permutation of the set $\{1, \ldots, \mu\}$ and subtract from
the initial point the elements of $\grob_{\prec_C}^{\sigma(i)}$,
for $i=1, \ldots, \mu$. Add this element to an empty list. After
each substraction by elements in $\grob_{\prec_C}^{\sigma(i)}$,
$i=1, \ldots, \mu$, remove from the list those elements dominated
by the new element. We prove in Section 3 that the result does not
depend on the permutation $\sigma$ chosen.

\setcounter{ej}{1}
\begin{ej}[Continuation]
\label{ej3} This example shows the abovementioned different ways
to compute paths from dominated solutions to any Pareto-optimal
solution. The vector $(9,4,9,3)$ is a feasible solution for
$MIP_{A,C}$ in the $(17,11)^t$-fiber. Figure \ref{fig2} shows the
sequence of Pareto-optimal points obtained from the feasible point
$(9,4,9,3)$ using the permutation $\sigma_1=(1,2)$ (on the left)
and using $\sigma_2=(2,1)$ (on the right).
\begin{figure}[h]
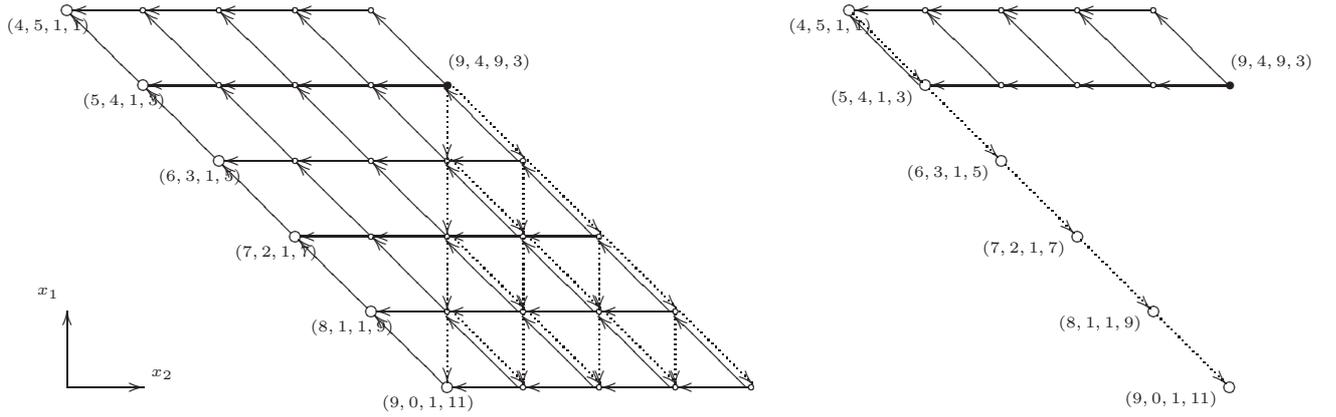

$$
\begin{array}{ll}
\xy
  \POS (40,0) *\cir<0.1pt>{}="a0-0"
, \POS (50,0) *\cir<0.1pt>{}="a1-0"
, \POS (40,10) *\cir<0.1pt>{}="a0-1"
, \POS (80,50) *\cir<1pt>{}="a1",
, \POS (70,50) *\cir<1pt>{}="a2"
, \POS (60,50) *\cir<1pt>{}="a3"
, \POS (50,50) *\cir<1pt>{}="a4"
, \POS (40,50) *\cir<2pt>{}="a5"
, \POS (90,40) *\dir{*}="a6"
, \POS (80,40) *\cir<1pt>{}="a7"
, \POS (70,40) *\cir<1pt>{}="a8"
, \POS (60,40) *\cir<1pt>{}="a9"
, \POS (50,40) *\cir<2pt>{}="a10"
, \POS (100,30) *\cir<1pt>{}="a11"
, \POS (90,30) *\cir<1pt>{}="a12"
, \POS (80,30) *\cir<1pt>{}="a13"
, \POS (70,30) *\cir<1pt>{}="a14"
, \POS (60,30) *\cir<2pt>{}="a15"
, \POS (110,20) *\cir<1pt>{}="a16"
, \POS (100,20) *\cir<1pt>{}="a17"
, \POS (90,20) *\cir<1pt>{}="a18"
, \POS (80,20) *\cir<1pt>{}="a19"
, \POS (70,20) *\cir<2pt>{}="a20"
, \POS (120,10) *\cir<1pt>{}="a21"
, \POS (110,10) *\cir<1pt>{}="a22"
, \POS (100,10) *\cir<1pt>{}="a23"
, \POS (90,10) *\cir<1pt>{}="a24"
, \POS (80,10) *\cir<2pt>{}="a25"
, \POS (130,0) *\cir<1pt>{}="a26"
, \POS (120,0) *\cir<1pt>{}="a27"
, \POS (110,0) *\cir<1pt>{}="a28"
, \POS (100,0) *\cir<1pt>{}="a29"
, \POS (90,0) *\cir<2pt>{}="a30"

\POS "a0-0" \ar^(1.25){{\text{\tiny{$x_1$}}}} "a0-1"
\POS "a0-0" \ar^(1.25){{\text{\tiny{$x_2$}}}} "a1-0"
\POS"a1" \ar "a2"
\POS"a2" \ar "a3"
\POS"a3" \ar "a4"
\POS"a4" \ar^(1.25){{\text{\tiny{$(4,5,1,1)$}}}} "a5"
\POS"a6" \ar "a7"
\POS"a7" \ar "a8"
\POS"a8" \ar "a9"
\POS"a9" \ar^(1.25){{\text{\tiny{$(5,4,1,3)$}}}} "a10"
\POS"a11" \ar "a12"
\POS"a12" \ar "a13"
\POS"a13" \ar "a14"
\POS"a14" \ar^(1.25){{\text{\tiny{$(6,3,1,5)$}}}} "a15"
\POS"a16" \ar "a17"
\POS"a17" \ar"a18"
\POS"a18" \ar "a19"
\POS"a19" \ar^(1.25){{\text{\tiny{$(7,2,1,7)$}}}} "a20"
\POS"a21" \ar "a22"
\POS"a22" \ar "a23"
\POS"a23" \ar "a24"
\POS"a24" \ar^(1.25){{\text{\tiny{$(8,1,1,9)$}}}}"a25"
\POS"a26" \ar "a27"
\POS"a27" \ar "a28"
\POS"a28" \ar "a29"
\POS"a29" \ar^(1.25){{\text{\tiny{$(9,0,1,11)$}}}} "a30"
\POS"a26" \ar@<0.3ex> "a21"
\POS"a21" \ar@<0.3ex> "a16"
\POS"a16" \ar@<0.3ex> "a11"
\POS"a11" \ar@<0.3ex> "a6"
\POS"a6" \ar_(0.1){{\text{\tiny{$(9,4,9,3)$}}}} "a1"
\POS"a21" \ar@<0.3ex>@{.>} "a26"
\POS"a16" \ar@<0.3ex>@{.>} "a21"
\POS"a11" \ar@<0.3ex>@{.>} "a16"
\POS"a6" \ar@<0.3ex>@{.>} "a11"
\POS"a7" \ar "a2"
\POS"a12" \ar "a7"
\POS"a17" \ar@<0.3ex> "a12"
\POS"a22" \ar@<0.3ex> "a17"
\POS"a27" \ar@<0.3ex> "a22"
\POS"a12" \ar@<0.3ex>@{.>} "a17"
\POS"a17" \ar@<0.3ex>@{.>} "a22"
\POS"a22" \ar@<0.3ex>@{.>} "a27"
\POS"a8" \ar "a3"
\POS"a13" \ar "a8"
\POS"a18" \ar "a13"
\POS"a23" \ar@<0.3ex> "a18"
\POS"a28" \ar@<0.3ex> "a23"
\POS"a18" \ar@<0.3ex>@{.>} "a23"
\POS"a23" \ar@<0.3ex>@{.>} "a28"
\POS"a9" \ar "a4"
\POS"a14" \ar "a9"
\POS"a19" \ar "a14"
\POS"a24" \ar "a19"
\POS"a29" \ar@<0.3ex> "a24"
\POS"a24" \ar@<0.3ex>@{.>} "a29"
\POS"a10" \ar "a5"
\POS"a15" \ar "a10"
\POS"a20" \ar "a15"
\POS"a25" \ar "a20"
\POS"a30" \ar "a25"
\POS"a6" \ar@{.>} "a12"
\POS"a12" \ar@{.>} "a18"
\POS"a18" \ar@{.>} "a24"
\POS"a24" \ar@{.>} "a30"
\POS"a11" \ar@{.>} "a17"
\POS"a17" \ar@{.>} "a23"
\POS"a23" \ar@{.>} "a29"
\POS"a16" \ar@{.>} "a22"
\POS"a22" \ar@{.>} "a28"
\POS"a21" \ar@{.>} "a27"

\endxy

&

\xy
\POS (80,50) *\cir<1pt>{}="b1"
, \POS (70,50) *\cir<1pt>{}="b2"
, \POS (60,50) *\cir<1pt>{}="b3"
, \POS (50,50) *\cir<1pt>{}="b4"
, \POS (40,50) *\cir<2pt>{}="b5"
, \POS (90,40) *\dir{*}="b6"
, \POS (80,40) *\cir<1pt>{}="b7"
, \POS (70,40) *\cir<1pt>{}="b8"
, \POS (60,40) *\cir<1pt>{}="b9"
, \POS (50,40) *\cir<2pt>{}="b10"
, \POS (60,30) *\cir<2pt>{}="b11"
, \POS (70,20) *\cir<2pt>{}="b12"
, \POS (80,10) *\cir<2pt>{}="b13"
, \POS (90,0) *\cir<2pt>{}="b14"

\POS"b1" \ar "b2"
\POS"b2" \ar "b3"
\POS"b3" \ar "b4"
\POS"b4" \ar^(1.25){{\text{\tiny{$(4,5,1,1)$}}}} "b5"
\POS"b6" \ar "b7"
\POS"b7" \ar "b8"
\POS"b8" \ar "b9"
\POS"b9" \ar "b10"
\POS"b6" \ar_(0.1){{\text{\tiny{$(9,4,9,3)$}}}} "b1"
\POS"b7" \ar "b2"
\POS"b8" \ar "b3"
\POS"b9" \ar "b4"
\POS"b10" \ar@<0.3ex> "b5"
\POS"b5" \ar@{.>}_(0.95){{\text{\tiny{$(5,4,1,3)$}}}} "b10"
\POS"b10" \ar@{.>}_(0.95){{\text{\tiny{$(6,3,1,5)$}}}} "b11"
\POS"b11" \ar@{.>}_(0.95){{\text{\tiny{$(7,2,1,7)$}}}} "b12"
\POS"b12" \ar@{.>}_(0.95){{\text{\tiny{$(8,1,1,9)$}}}} "b13"
\POS"b13" \ar@{.>}_(0.95){{\text{\tiny{$(9,0,1,11)$}}}} "b14"\
\endxy
\end{array}
$$
\caption{\label{fig2}Two different ways to compute paths from
$(9,4,9,3)$ to the Pareto-optimal solutions in its fiber.}
\end{figure}

\end{ej}
\begin{remark}
With the partial order, $\prec_C$, induced by $C$ a directed path
from a dominated point $\alpha$ to each Pareto-optimal point
$\beta$ in a fiber, applying the above method, cannot pass through
any lattice point in this fiber more than $\mu$ times (recall that
$\mu$ is the number of maximal chains in $\grob_{\prec_C}$). This
implies that obtaining the Pareto-optimal solutions of a given
$MIP_{A,C}$ using $\grob_{\prec_C}$ cannot cycle.
\end{remark}

\section{Test families and Partial Gr\"obner bases}
\label{sect2} In the previous section we motivate the importance
of having a test family for $MIP_{A,C}$ since this structure
allows us obtaining the entire set of Pareto-optimal solutions of
the above family of multiobjective integer programs (when the
right-hand side varies). Our goal in this section is to provide
the necessary tools to construct test families for any
multiobjective integer problem. Our construction builds upon an
extension of Gröbner bases on partial orders.

In order to introduce this structure we define the reduction of a
pair $(g,h) \in \Z^n \times \Z^n_+$ by a finite set of ordered
pairs in $\Z^n \times \Z^n_+$.  Given is a collection $\grob_C
\subseteq \Z^n \times \Z^n_+$ where $\grob_C = \{ (g_{1}, h_{1}),
\ldots, (g_{l}, h_{l}) : h_{k+1} \prec_C h_{k}, k = 1, \ldots,
l-1\}$.

The reduction of $(g,h)$ by $\grob_C$ consists of the process
described in Algorithm \ref{preduction}.

\begin{algorithm}[H]
\label{preduction}

\SetLine \SetKwInOut{Input}{input} \SetKwInOut{Output}{output}

\Input{$R=\{(g,h)\}$, $S=\{(g,h)\}$} For each $(\tilde{g}, \tilde{h}) \in S:$
\Repeat{$\{i : \tilde{h} - h_i \geq 0\} = \emptyset$}{
\eIf{$\tilde{h}-g_i$ and $\tilde{h}-\tilde{g}$ are comparable by
$\prec_C$}{$R_o = \{(\tilde{g}-g_i, \max_{\prec_C}
\{\tilde{h}-\tilde{g}_i,\tilde{h}-\tilde{g}\})\}$}{$R_o
=\{(\tilde{g}-g_i, \tilde{h}-g_i) , (\tilde{g}-g_i
,\tilde{h}-\tilde{g})\}$}

For each $r \in R_o$ and $s \in R$:

\If{$r \prec_C s$}{ $R:= R \backslash \{s\}$;}

$S:=R_o$

 $R := R \cup R_o$;
} \Output{$R$, the partial reduction set of $(g,h)$ by $\grob_C$}
 \caption{Partial reduction algorithm}
\end{algorithm}

The above reduction process extends to the case of a finite
collection of ordered sets of pairs in $\Z^n \times \Z^n_+$ by
establishing the sequence in which the sets of pairs are
considered. We denote by $pRem((g,h),(\grob_i))_\sigma$ the
reduction of the pair $(g,h)$ by the family $\{\grob_i\}_{i=1}^t$
for a fixed sequence of indices $\sigma$. The following result
allows us to consider any sequence of indices for this process,
since it establishes that the partial reduction does not depend on
the chosen sequence.

\begin{theo} Let $\grob$ be a finite set in
$\Z^n \times \Z^n_+$, whose maximal chains are $\grob_1, \ldots,
\grob_t$, and $\sigma$, $\sigma^{'}$ two permutations of the
indices $(1, \ldots, t)$. Then,
$$
pRem((g,h),\grob)_\sigma = pRem((g,h),\grob)_{\sigma^{'}}
$$
for each $(g, h) \in \Z^n \times \Z^n_+$.
\end{theo}
\begin{proof}

Consider $\Lambda_\sigma:= \{ \tilde{g} : \tilde{g} = g -
\dsum_{i=1}^t \dsum_{j=1}^{k_{\sigma(i)}}
\lambda_\sigma(i)\,g^{\sigma(i)}_j\}$, where $\grob_i = \{(g_j^i,
h_i^j): j=1\ldots, k_i\}$. It is clear that the elements in
$\Lambda_\sigma$ does not depend on the permutation $\sigma$,
since reordering the sums does not give new elements.

The elements in $pRem((g,h),\grob)_\sigma$ are the element in
$\Lambda_\sigma$ deleting the comparable largest ones. Then, since
$\Lambda_\sigma = \Lambda_{\sigma^{'}}$, $pRem((g,h),\grob)_\sigma
= pRem((g,h),\grob)_{\sigma^{'}}$.

\end{proof}

From now on, we denote by $pRem((g,h), \grob)$ the set of
remainders of $(g,h)$ by the family $\grob = \{\grob_i\}_{i=1}^t$
for any sequence of indices.

The reduction of a pair that represents a feasible solution, by a
test family, gives the entire set of Pareto-optimal solutions. In
order to obtain that test family, we introduce the notion of
p-Gr\"obner basis. This concept has been motivated by the fact
that in the case where the ordering induced in $\N^n$ by a single
cost vector is total, a Gr\"obner basis is a test set for the
family of integer programs $IP_{A,c}$ (see \cite{conti-traverso91}
or \cite{thomas95} for extended details). In the single objective
case the Buchberger algorithm computes the Gr\"obner basis.
However, in the multiobjective case the cost matrix induces a
partial order, so division or the Buchberger algorithm are not
applicable. Using the above reduction algorithm we present in this
section an ``à la'' Buchberger algorithm to compute the so called
p-Gr\"obner basis to solve MOILP problems.
\begin{defi}[Partial Gr\"obner basis]
A family $\grob = \{\grob_1, \ldots, \grob_t\} \subseteq I_A$ is a
partial Gr\"obner basis (p-Gr\"obner basis) for the family of
problems $MIP_{A,C}$, if $\grob_1, \ldots, \grob_t$ are the
maximal chains for the partially ordered set $\dbigcup_{i=1}^t
\grob_i$ and for any $(g,h) \in \Z^n\times\Z^n_+$:
$$
g \in Ker(A) \Longleftrightarrow pRem((g,h),\grob)= \{0\}.
$$

A p-Gr\"obner basis is said to be \textit{reduced} if every
element at each maximal chain cannot be obtained by reducing any
other element of the same chain.
\end{defi}

Given a p-Gr\"obner basis, computing a reduced p-Gr\"obner basis
is done by deleting the elements that can be reduced by other
elements in the basis. After the removing process, the family is a
p-Gr\"obner basis having only non redundant elements. It is easy
to see that the reduced p-Gr\"obner basis for $MIP_{A,C}$ is
unique and minimal, in the sense that no element can be removed
from it maintaining the p-Gr\"obner basis structure.

This definition clearly extends to p-Gr\"obner bases for the ideal
$I_A$ induced by $A$, once we fix the partial order, $\prec_C$,
induced by $C$.

The goal of this paper is to present algorithms to solve
multiobjective problems analogous to the methods that solve the
single objective case, using usual Gr\"obner basis. These methods
are based on computing the reduction of a feasible solution by the
basis. The key for that result is the fact that the reduction of
any pair of feasible solutions is the same, therefore the
algorithm is valid for any initial feasible solution. The
following lemma assures the same statement for the multiobjective
case and p-Gr\"obner bases.

\begin{lemma}\label{lemma1}
Let $\grob$ be the reduced p-Gr\"obner basis for $MIP_{A,C}$ and
$\alpha_1, \alpha_2$ two different feasible solutions in the same
fiber of $MIP_{A,C}$. Then, $pRem((\alpha_1,\alpha_1),\grob) =
pRem((\alpha_2,\alpha_2),\grob)$.
\end{lemma}
\begin{proof}
Let $(\beta, \beta) \in pRem((\alpha_1,\alpha_1),\grob)$, then
$\beta - \alpha_2$ is in the same fiber and it cannot be reduced,
so $(\beta,\beta) \in pRem((\alpha_2,\alpha_2),\grob)$.
\end{proof}
The following theorem states the relationship between the three
structures introduced before: test families, reduced p-Gr\"obner
bases and the family $\grob_{\prec_C}$.
\begin{theo}
\label{testfam-pgrobner} The reduced p-Gr\"obner basis for
$MIP_{A,C}$ is the unique minimal test family for $MIP_{A,C}$.
Moreover, $\grob_{\prec_C}$, introduced in \eqref{gc}, is the
reduced p-Gr\"obner basis for $MIP_{A,C}$.
\end{theo}

\begin{proof}
Let $\grob = \{\grob_1, \ldots, \grob_t\}$ be the reduced
p-Gröbner basis for $MIP_{A,C}$. By definition of p-Gröbner basis,
it is clear that each $\grob_i$ is totally ordered by its second
component with respect to $\prec_C$ (Condition 1). Condition 2
follows because for each $i$ and for each $(g,h) \in \grob_i
\subseteq \Z^n \times \Z^n_+$, clearly $pRem((g,h), \grob) =
\{0\}$, so $g \in Ker(A)$ and then $A(h-g) = Ah$.

Now, let $x\in \N^n$ be a dominated solution for $MIP_{A,C}(b)$
then there is a Pareto-optimal solution, $\beta$, such that $\beta
\prec_C x$. By Lemma \ref{lemma1}, $pRem((x,x), \grob) =
pRem((\beta, \beta), \grob)$, and by construction of the set of
partial remainders, $\beta \in pRem((\beta, \beta), \grob)$, and
then $x \not\in pRem((x,x),\grob)$. This implies that $(g,h) \in
\grob_i$ must exist such that $x - g_i \prec_C x$, for some $i \in
\{1, \ldots, t\}$.

On the other hand, if $x$ is a Pareto-optimal solution for
$MIP_{A,C}(b)$, $x \in pRem((x,x),\grob)$, and then, there exists
no $(g,h)$ in any $\grob_i$ such that $x-g \prec_C x$. Therefore,
for every $i$ and for each $(g,h) \in \grob_i$, either $x-g$ is
infeasible or incomparable with $x$.

Minimality is due to the fact that removing an element from the
reduced p-Gröbner basis, that is the minimal partial Gröbner basis
that can be built for $MIP_{A,C}$ we cannot guarantee to have a
test family because it may exist a pair $(g,h) \in \Z^n\times
\Z^n_+$ with $g \in Ker(A)$ that cannot be reduced to the zero
set.

The second statement of the theorem follows from Corollary
\ref{grobC}.
\end{proof}

In the following we describe an extended algorithm to compute a
p-Gr\"obner basis for $I_A$, with respect to the partial order
induced by $C$. First, for $(g,h), (g^{'}, h^{'})$ in $\Z^n \times
\Z_+^n$ we denote by $S^1((g,h),(g^{'},h^{'}))$ and
$S^2((g,h),(g^{'},h^{'}))$ the pairs
$$ S^1((g,h),(g^{'},h^{'})) =
\left\{ \begin{array}{ll}
(g - g^{'} - 2(h-h^{'}), \gamma+g-2h) & \mbox{if $\gamma+g-2h \prec_C \gamma+g^{'}-2h^{'}$} \\
(g^{'} - g - 2(h^{'}-h),\gamma+g^{'}-2h^{'}) & \mbox{if $\gamma+g^{'}-2h^{'} \prec_C \gamma+g-2h$}\\
(g - g^{'} - 2(h-h^{'}),\gamma+g-2h) & \mbox{if
$\gamma+g^{'}-2h^{'}$ and $\gamma+g-2h$ are incomparable}
\end{array}\right.
$$
$$ S^2((g,h),(g^{'},h^{'})) =
\left\{ \begin{array}{ll}
(g - g^{'} - 2(h-h^{'}), \gamma+g-2h) & \mbox{if $\gamma+g-2h \prec_C \gamma+g^{'}-2h^{'}$} \\
(g^{'} - g - 2(h^{'}-h),\gamma+g^{'}-2h^{'}) & \mbox{if $\gamma+g^{'}-2h^{'} \prec_C \gamma+g-2h$}\\
(g^{'} - g - 2(h^{'}-h),\gamma+g^{'}-2h^{'}) & \mbox{if
$\gamma+g^{'}-2h^{'}$ and $\gamma+g-2h$ are incomparable}
\end{array}\right.
$$
where $\gamma \in \N^n$ and $\gamma_i = \max\{h_i, h^{'}_i\}$,
$i=1, \ldots, n$.

The pairs $S^1((g,h),(g^{'},h^{'}))$ and
$S^2((g,h),(g^{'},h^{'}))$ are called $1-Svector$ and $2-Svector$
of $(g,h)$ and $(g^{'},h^{'})$, respectively.

The notation is due to the analogy with the algebraic-geometrical
notion of S-polynomial for a pair of polynomials with a given term
order. Since we consider a partial order, it may happen that in
the standard construction of a Svector (see \cite{thomas95}), we
cannot decide which is the leading term. Therefore, in our
definitions of Svectors we follow the standard construction but we
must consider all possible combinations of leading terms, with
respect to the partial order $\prec_C$. The following lemma is
used in the proof of our extended criterion and it is an
adaptation of the analogous result for total orders and usual
S-polynomials.

In the following, we denote by $leadmon_C(f)$ the set of leading
monomials with respect to the order induced by $C$, for any
multivariate polynomial $f\in K[x_1, \ldots, x_n]$ and by
$1$-Spolynomial and $2$-Spolynomial the binomial transcriptions of
$1$-Svector and $2$-Svector (recall the equivalence between the
pairs $(u,v)$ and the binomial $x^{u-v}-x^u$ if $u$ is dominated
by $v$).

\begin{lemma}\label{lemma:pspol}
Let $f_1, \ldots, f_s$ $\in$ $K[x_1, \ldots, x_n]$ be such that
there exists $p \in \dbigcap_{i=1}^s leadmon_C(f_i)$. Let $f=
\dsum_{i=1}^s c_i\,f_i$ with $c_i \in K$. If there exists $q \in
leadmon_C(f)$ such that $q \prec_C p$, then $f$ is a linear
combination with coefficients in $K$ of the $k$-Spolynomial,
$k=1,2$,  of $f_i$ and $f_j$, $1 \leq i < j \leq s$.
\end{lemma}
\begin{proof}
By hypothesis, $f_i= a_i\,p + \mbox{ \textit{other smaller or
incomparable terms}}$, with $a_i \in K$, for all $i$. Then, $f$
can be rewritten as $f= \dsum_{i=1}^s c_i\,f_i = \dsum_{i=1}^s
c_i\,a_i\, p  + \mbox{ \textit{other smaller or incomparable
terms}}$. Since $q \prec_C p$, then $\dsum_{i=1}^s c_i\,a_i = 0$.

By definition, for $k=1,2$, $S^k((f_i,p), (f_j,p)) =
\frac{1}{a_i}\,f_i - \frac{1}{a_j}\,f_j$, thus,
\begin{eqnarray*}
\begin{split}
f &= c_1\,f_1 + \cdots + c_s\,f_s \\
  &= c_1\,a_1 (\frac{1}{a_1}\,f_1) + \cdots + c_s\,a_s
  (\frac{1}{a_s}\,f_s)\\
  &= c_1\,a_1 (\frac{1}{a_1}\,f_1 - \frac{1}{a_2}\,f_2) +  (c_1\,a_1 + c_2\,a_2)\,(\frac{1}{a_2}\,f_2 - \frac{1}{a_3}\,f_3)
   + \cdots\\
&+ (c_1\,a_1 + \cdots +
c_{s-1}\,a_{s-1})\,(\frac{1}{a_{s-1}}\,f_{s-1} -
  \frac{1}{a_s}\,f_s) + (c_1\,a_1 + \cdots + c_{s}\,a_{s})\,
  \frac{1}{a_s}\,f_s \\
  &= d_1\,S((f_1,p),(f_2,p)) + \cdots + d_{s-1}\,S((f_{s-1},p),(f_s,p)) + \dsum_{i=1}^s
  c_i\,a_i = \dsum_{i=1}^{s-1} d_i\,S((f_{i},p),(f_{i-1},p)).
\end{split}
\end{eqnarray*}
where $d_i = \dsum_{j=1}^i\,c_j\,a_j$. This proves the lemma.
\end{proof}

The algorithm to compute standard Gr\"obner bases is based on the
Buchberger criterion, whose analogous for a partial order is the
following.

\begin{theo}[Extended Buchberger's criterion]
\label{th:extendedbuch} Let $\grob=\{\grob_1, \ldots, \grob_t\}$
with $\grob_i \subseteq I_A$ for all $i=1, \ldots, t$, be the
maximal chains of the partially ordered set $\{g_i : g_i \in
\grob_i, \text{ for some } i=1, \ldots, t\}$. Then the following
statements are equivalent:
\begin{enumerate}
\item $\grob$ is a p-Gr\"obner basis for the family $MIP_{A,C}$.
\item For each $i, j = 1, \ldots, t$ and $(g,h)\in \grob_i$, $(g^{'},h^{'})\in \grob_j$, $pRem(S^k((g,h),(g^{'},h^{'})), \grob)
= \{0\}$ , for $k=,1,2$.
\end{enumerate}
\end{theo}
\begin{proof} The original Buchberger
criterion was stated in a polynomial language. Therefore, we adapt
our notation to follow the line of that proof. Each pair $\{u,v\}$
is identified with the binomial $x^u-x^v$, in the polynomial ring
$\Z[x_1,\ldots,x_n]$, and our set $I_A$, with $\Im_A=\langle x^u -
x^v : u-v \in Ker(A) \rangle$. The definition of partial
remainders, $pRem$, is adapted accordingly. With these changes in
the notation, the set $setlm(\{u,v\})$ is identified with the
elements in $leadmon_C(x^u - x^v)$.

Let $\grob$ be a p-Gr\"obner basis for $I_A$, $i, j \in \{1,
\ldots, t\}$ and $(g,h)\in \grob_i$, $(g^{'},h^{'})\in \grob_j$.
Then, $S^k((g,h),(g^{'},h^{'}))$, for $k=1, 2$, is in $I_A$, so by
definition of p-Gr\"obner basis,
$pRem(S^k((g,h),(g^{'},h^{'})),\grob) = \{0\}$.

Conversely, assume that for each $(\widetilde{g},\widetilde{h})\in
\grob_i$ and $(g^{'},h^{'})\in \grob_j$,
$pRem(S^k((\widetilde{g},\widetilde{h}),(g^{'},h^{'})), \grob)
\neq \{0\}$ , for $k=,1,2$. Let $(g,h) \in \Z^n\times\Z^n_+$ with
$g \in Ker(A)$. We define $f = x^{h}-x^{g-h} \in \Z[x_1, \ldots,
x_n]$, and we denote by $\grob^*$ the polynomial set associated
with $\grob$.

Then, $f$ can  be written as a linear combinations of all the
elements in $\grob^*$ (this representation is not unique):
    $$
    f= \dsum_{i=1}^m h_i\,g_i.
    $$
Let $X = \{ X_1, \ldots, X_N\}$ be the set of maximal elements of
the set $\{H_i\,G_i : H_i \in leadmon_C(h_i), G_i \in
leadmon_C(g_i)\}$, with respect to the partial order $\prec_C$.

If $X= leadmon_C(f)$, the polynomial $f$ can be partially reduced
by the elements in $\grob$. This proves the result.

Otherwise, assume that $l \in leadmon_C(f)\backslash X$. Then, $l$
comes from some simplification of the linear combination defining
$f$. Then, the construction ensures that it must exist at least
one element, $X_i \in X$, such that $l \prec_C X_i$.

Set $S= \{ j : H_j\,G_j = X_i \mbox{ with $H_i \in leadmon_C(h_i),
G_i \in leadmon_C(g_i)$}\}$. For any $j \in S$, we write $h_j =
H_j + \mbox{ \textit{other terms} }$ and $g= \dsum_{j\in S}
H_j\,g_j$. Then, $X_i \in leadmon_C(H_j\,g_j)$, for all $j \in S$.
However, by hypothesis there exists $G \in leadmon_C(g)$, with $G
\prec X_i$.

Hence, by Lemma \ref{lemma:pspol}, there exists $d^k_{s,r} \in K$
such that:
 $$
    g= \dsum_{k=1}^2 \dsum_{r,s\in S, r\neq s, g_s,g_r \in \grob_k}
    d_{s,r}\,S^k(X_s\,g_s,X_r\,g_r).
$$

Now, for any $r$, $s\in S$, $X_i = lcm(L_r, L_s)$ for some $L_r
\in leadmon_C(H_r\,g_r)$ and $L_s \in leadmon_C(H_s\,g_s)$, so:

    \begin{eqnarray*}
    \begin{split}
    S^k((H_r\,g_r, L_r), (H_s\,g_s, L_s)) &=  \dfrac{X_i}{L_r} \,H_r\,g_r
    - \dfrac{X_i}{L_s} \,H_s\,g_s\\
    & = \dfrac{X_i}{l_r} \,g_r
    - \dfrac{X_i}{l_s}\,g_s = \dfrac{X_i}{H_{r,s}} \,S^k((g_r,l_r),(
    g_s,l_s)
    \end{split}
    \end{eqnarray*}
where $l_r = \frac{l_r}{H_r}$, $l_s = \frac{l_r}{H_s}$ and
$H_{r,s} := lcm(lp(g_r),lp(g_s))$.

By hypothesis, $pRem(S^k(g_r,g_s), \grob) = \{0\}$. From the last
    equation we deduce that:
$$
pRem(S^k(H_r\,g_r,H_s\,g_s), \grob)= \{0\}
$$
this gives a representation:
$$
S^k(H_r\,g_r,H_s\,g_s) = \dsum_{\nu}\,h_{r,s}^\nu\,g_\nu
$$
with $g_\nu \in \grob$:
$$
\dmax_{\nu} \{H_{r,s}^\nu \,G^\nu \mbox{: $H_{r,s}^\nu \in
leadmon_C(h_{r,s}^\nu), G^\nu \in
    leadmon_C(g_\nu)$}\} =
    leadmon_C(S(H_r\,g_r,H_s\,g_s)^{k})=:S_{r,s}^k.
$$

By construction of S-polynomials, we have that there exists $p \in
S_{r,s}^k$ such that $p \prec_C X_i$, so, substituting these
expressions into $g$ above and using that $f= \dsum_{j\not\in S}
h_j\,g_j + \dsum_{j \in S} h_j\,g_j = \dsum_{j\not\in S} h_j\,g_j
+ g = \dsum_{j\not\in S} h_j\,g_j + \dsum_{r,s}
d_{r,s}\,S(H_s\,g_s, H_r\,g_r) = \dsum_{j\not\in S} h_j\,g_j +
\dsum_{r,s}\dsum_{\nu} h_{r,s}^\nu\,g_\nu$, we have expressed $f$
as:
$$
f= \dsum_i \,h_i^{'}\,g_i
$$
with one leading term, $p$, smaller than $X_i$. However, this is a
contradiction proving the theorem.
\end{proof}

This criterion (the one in Theorem \ref{th:extendedbuch}) allows
us to describe a geometric algorithm which constructs a
p-Gr\"obner basis $\grob_C$ for $MIP_{A,C}$, and then a test
family for that family of multiobjective problems.

The first approach to compute a p-Gr\"obner basis for a family of
multiobjective programs, is an algorithm based on Conti and
Traverso method for the single objective case
\cite{conti-traverso91}. For this algorithm, the key is
transforming the given multiobjective program into another one
where computation is easier and an initial set of generators for
$I_A$  are known.

Notice that finding an initial set of generators for $I_A$ can be
done by a straightforward modification of the Big-M method (see
details, e.g. in \cite{bazaraa93}).

Given the program $MIP_{A,C}(b)$, we consider the associate
extended multiobjective program, $EMIP_{A,C}(b)$ as the problem
$MIP_{\widetilde{A},\widetilde{C}}(b)$ where $\widetilde{A} = \
\left(\begin{array}{c|c|c}
 & -1 &\\
Id_{m}& \vdots & A\\
 & -1 &
\end{array}\right) \in \Z^{m\times (m+1+n)}$, $\widetilde{C} =
(M\cdot\mathbf{1} |C) \in \Z^{(m+1+n)\times k}$, $Id_m$ stands for
the $m\times m$ identity matrix, $M$ is a large constant and
$\mathbf{1}$ is the $(m+1)\times k$ matrix whose components are
all $1$. This problem adds $m+1$ new variables, whose weights in
the multiobjective function are big, and so, solving this extended
minimization program allows us to solve directly the initial
program $MIP_{A,C}$. Indeed, any feasible solution to the original
problem is a feasible solution to the extended problem with the
first $m$ components equal to zero, so any feasible solution of
the form $(0, \stackrel{m+1}{\ldots}, 0, \alpha_1, \ldots,
\alpha_n)$ is non-dominated, upon the order
$\prec_{\widetilde{C}}$, by any solution without zeros in the
first $m$ components. Then, computing a p-Gr\"obner basis for the
extended program, allows us detecting infeasibility of the
original problem. Furthermore, a trivial feasible solution,
$\widetilde{\mathbf{x}}_0 = (b_1, \ldots, b_m, 0,
\stackrel{n+1}{\ldots}, 0)$, is known and the initial set of
generators for $I_A$ are given by $\{\{M_i-P_i, M_i\}: i=0\ldots,
n\}$ where $M_i = (a_{1i}-\min\{0, \min_j\{a_{ji}\}\}, \ldots,
a_{mi}-\min\{0, \min_j\{a_{ji}\}\},-\min\{0, \min_j\{a_{ji}\}\},
0, \stackrel{n}{\ldots}, 0)$, $P_i = (0, \stackrel{m+1}{\ldots}, 0
| e_i)$, for all $i = 1, \ldots, n$, $M_0 =(1,
\stackrel{m+1}{\ldots}, 1, 0, \stackrel{n}{\ldots},0)$ and $P_0 =
\mathbf{0}$, $M_i, P_i, M_0, P_0 \in \Z_+^{n+m+1}$ (see
\cite{adams94} for further details).

\begin{algorithm}[H]
\label{pbuchberger1} \SetLine

\SetKwInOut{Input}{input} \SetKwInOut{Output}{output}

\Input{$F_1=\{M_0, M_1, \ldots, M_n\}$ and $F_2=\{P_0, P_1, \ldots, P_n\}$,
$M_i = (a_{1i}-\min\{0, \min_j\{a_{ji}\}\}, \ldots, a_{mi}-\min\{0, \min_j\{a_{ji}\}\},-\min\{0, \min_j\{a_{ji}\}\}, 0, \stackrel{n}{\ldots}, 0)$ ($i>0$)\\
$P_i = (0, \stackrel{m+1}{\ldots}, 0 | e_i) \in \N^{m+n+1}$ ($i>0$)\\
$M_0 =(1, \stackrel{m+1}{\ldots}, 1, 0, \stackrel{n}{\ldots},0)$\\
$P_0 = (0, \stackrel{n+m+1}{\ldots}, 0)$.}

\Repeat{$R^k = \{0\}$ for every pairs}{ Compute, $\grob_1, \ldots,
\grob_t$, the maximal chains for $\grob = \phi(\mathcal{F}(F_1,
F_2))$.

\For{$i,j \in \{1, \ldots, t\}$, $i\neq j$, and each pair $(g,h)
\in \grob_i$, $(g',h') \in \grob_j$}{Compute $R^k =
pRem(S^k((g,h),(g^{'},h^{'})),\grob)$, $k=1,2$.

\eIf{$R^k = \{0\}$}{Continue with other pair.}{Add
$\phi(\mathcal{F}(r))$ to $\grob$, for each $r\in R^k$.}

}
 }
\Output{$\grob = \{\grob_1, \ldots, \grob_Q\}$} p-Gröbner basis
for $I_A$ with respect to $\prec_C$.
 \caption{Partial Buchberger algorithm I}
\end{algorithm}

Then, we can state the following result.

\begin{theo}
Let  $\grob = \{\grob_i\}_{i=1}^t$ be a p-Gr\"obner basis for
$EMIP_{A,C}$. If $(0, \stackrel{m+1}{\ldots}, 0, \alpha_1, \ldots,
\alpha_n) \in pRem((0, \stackrel{m+1}{\ldots}, 0, b_1, \ldots,
b_n), \grob)$, then $\alpha=(\alpha_1, \ldots, \alpha_n)$ is a
Pareto-optimal solution for $MIP_{A,C}(b)$. The entire set of
Pareto-optimal solutions of $MIP_{A,C}(b)$ can be computed using
the above construction. Moreover, if there are no $\alpha$ in the
set $pRem((0, b),\grob)$ whose $m+1$ first components are zero
$MIP_{A,C}(b)$ is infeasible.
\end{theo}
\begin{proof}
Let $\alpha$ be a vector obtained by successive reductions over
$\grob$.
 It is clear that $\alpha$ is feasible because $(\mathbf{0},\alpha)$ is in
  the set of remainders of $(\mathbf{0}, \beta)$ and then, in the same
  fiber. Besides, $\alpha$ is a Pareto-optimal
  solution because $\grob$ is a test family for the problem
  (Theorem \ref{testfam-pgrobner}).

Now, if $\beta^*$ is a Pareto-optimal solution, by Lemma
\ref{lemma1} $pRem((\beta^*, \beta^*),\grob)) = pRem((\beta,
\beta),\grob))$, but since $\beta^*$ is a Pareto-optimal solution,
it cannot be reduced so $(\beta^*,\beta^*) \in pRem((\beta^*,
\beta^*),\grob))$, and then, also to the list of partial
remainders of $(\beta, \beta)$ by $\grob$.
\end{proof}

Hosten and Sturmfels \cite{hosten-sturmfels95} improved the method
by Conti and Traverso to solve single-objective programs using
standard Gr\"obner bases. Their improvement comes from the fact
that it is not necessary to increase the number of variables in
the problem, as Conti and Traverso's algorithm does.
Hosten-Sturmfels's algorithm allows decreasing the number of steps
in the computation of the Gr\"obner basis, but on the other hand,
it needs an algorithm to compute an initial feasible solution,
that in Conti and Traverso algorithm was trivial. We have modified
this alternative algorithm to be used to compute the entire set of
Pareto-optimal solutions. The first step in the algorithm is
computing an initial basis for the polynomial toric ideal
$\Im_A=\langle x^u - x^v : u-v \in Ker(A)\rangle$, that we are
identifying with $I_A$. This step does not depend on the order
induced by the objective function, so it can be used to solve
multiobjective problems. Details can be seen in
\cite{hosten-sturmfels95}. Algorithm \ref{systemofgenerators}
implements the computation of the set of generators of $\Im_A$.
This procedure uses the notion of LLL-reduced basis (see
\cite{lll82} for further details). In addition, we use a
$\omega$-graded reverse lexicographic term order,
$\prec^{gr_i}_{\omega}$, induced by $x_{i+1} > \cdots
> x_{i-1} > x_i$ (with $x_{n+1} := x_1$), that is defined as follows:
$$
\alpha \prec^{gr_i}_{\omega} \beta : \Longleftrightarrow \left\{
\begin{array}{ll}
\sum_{i=1}^n \omega_i \alpha_i < \sum_{i=1}^n \omega_i \beta_i &
\mbox{or}\\
\sum_{i=1}^n \omega_i \alpha_i = \sum_{i=1}^n \omega_i \beta_i &
\mbox{and $\alpha \prec_{lex} \beta$}
\end{array}\right.
$$
where $\omega \in \R_+^n$ is chosen such that $x_{i+1} > \cdots
> x_{i-1} > x_i$.

\begin{algorithm}[H]
\label{systemofgenerators} \SetLine

\SetKwInOut{Input}{input} \SetKwInOut{Output}{output}

\Input{$A \in \Z^{m\times n}$}

\begin{enumerate}

\item Find a lattice basis $\mathcal{B}$ for $Ker(A)$ (using the
Hermite Normal Form).
\item Replace $\mathcal{B}$ by the LLL-reduced lattice basis
$\mathcal{B}_{red}$ in the sense of L\`ovasz (see \cite{lll82} for more details).\\
Let $J_0 := \langle x^{u_+} - x^{u_-} : u \in \mathcal{B}_{red}
\rangle$.

\For{$i=1, \ldots, n$}{Compute $J_i =(J_{i-1}: x_i^\infty)$ as:
\begin{enumerate}
    \item Compute $\grob_{i-1}$ the reduced
    Gr\"obner basis for $J_{i-1}$ with respect to $\prec^{gr_i}_{\omega}$.
    \item Divide each element $f \in \grob_{i-1}$ by the highest
    power of $x_i$ that divides $f$.
\end{enumerate}}

\end{enumerate}
\Output{$\Im_A := J_n = \{x^{u_1}-x^{v_1}, \ldots,
x^{u_s}-x^{v_s}\}$ system of generators for $I_A$.}
 \caption{\texttt{setofgenerators($A$)}}
\end{algorithm}

$\Im_A$ consists of binomials $x^{u_i} - x^{v_i}$ with $u_i-v_i\in
Ker(A)$, for $i=1, \ldots, s$. Coming back to our notation, each
binomial, $x^u - x^v$, in $\Im_A$ is identified with $\{u,v\} \in
I_A$, so computing a set of generators for $\Im_A$ gives us, in
some sense, a finite number of generators for the set that
represents the constraints matrix. We compute in the next step a
partial Gr\"obner basis from the initial sets $F_1=\{u_1, \ldots,
u_s\}$ and $F_2=\{v_1, \ldots, v_s\}$ using our extended
Buchberger algorithm:

\begin{algorithm}[H]
\label{pbuchberger2} \SetLine

\SetKwInOut{Input}{input} \SetKwInOut{Output}{output}

\Input{$F_1=\{M_1, \ldots, M_r\}$ and $F_2=\{P_1, \ldots, P_r\}$.}

\Repeat{$R^k = \{0\}$ for every pairs}{ Compute, $\grob_1, \ldots,
\grob_t$, the maximal chains for $\grob = \phi(\mathcal{F}(F_1,
F_2))$.

\For{$i,j \in \{1, \ldots, t\}$, $i\neq j$, and each pair $(g,h)
\in \grob_i$, $(g',h') \in \grob_j$}{Compute $R^k =
pRem(S^k((g,h),(g^{'},h^{'})),\grob)$, $k=1,2$.

\eIf{$R^k = \{0\}$}{Continue with other pair.}{Add
$\phi(\mathcal{F}(r))$ to $\grob$, for each $r\in R^k$.}

}
 }
\Output{$\grob = \{\grob_1, \ldots, \grob_Q\}$} p-Gröbner basis
for the set spanned by $\{\{M_i, P_i\}: i=1,\ldots, r\}$ with
respect to $\prec_C$.
 \caption{\texttt{pgrobner($F_1,F_2$)}}
\end{algorithm}

Once we have obtained  the partial Gröbner basis using the above
algorithm, we can compute the entire set of Pareto-optimal
solutions for $MIP_{A,C}(b)$ by the following algorithm:

\begin{algorithm}[H]
\label{pos} \SetLine

\SetKwInOut{Input}{input} \SetKwInOut{Output}{output}

\Input{$MIP_{A,C}(b)$}
\begin{description}
\item[\textsc{Step 1. }] Compute an initial feasible solution, $\alpha_o$, for $MIP_{A,C}(b)$.
It consists of finding a solution for the diophantine system of
equations $Ax = b$, $x \in \Z^n$.
\item[\textsc{Step 2. }] Compute a system of generators for $I_A$:
$\{\{u_i,v_i\}: i=1, \ldots, s\}$, using
\verb"setofgenerators("$A$\verb")".
\item[\textsc{Step 3. }] Compute the partial reduced Gr\"obner basis for $MIP_{A,C}$,
$\grob_C = \{\grob_1, \ldots, \grob_t\}$, using \verb"pgrobner("
$F_1,F_2$ \verb")", where $F_1 = \{u_i : i=1, \ldots, r\}$ and
$F_2 = \{v_i : i=1, \ldots, r\}$.
\item[\textsc{Step 4. }] Calculate the set of partial remainders: $R:=
pRem(\alpha_o,\grob_C)$.
\end{description}
\Output{Pareto-optimal Solutions : $R$.}
 \caption{Pareto-optimal solutions computation for $MIP_{A,C}(b)$}
\end{algorithm}

There are some interesting cases where our methodology is highly
simplified due to the structure of the set of constraints. One of
these cases is when the dimension of the set of constraints is
$n-1$. The next remark explains how the algorithm simplifies in
this case.

\begin{remark}
\label{remark:dimkernel} Let $A$ be a $m \times n$ integer matrix
with rank $n-1$. Then, since $dim(Ker(A)) =1$, the system of
generators for $I_A$ $(${\bf Step 2}$)$ has just one element,
$(g,h)$, and the p-Gr\"obner basis $(${\bf Step 3}$)$ is the
family $\grob = \{ \{(g,h)\}\}$ because no Svector appears during
the computation of the Buchberger algorithm. In this case,
Pareto-optimal solutions are obtained as partial remainders of an
initial feasible solution $(\alpha, \alpha)$ by $(g,h)$, i.e., the
entire set of Pareto-optimal solutions is a subset of $\Gamma = \{
\alpha - \lambda g : \lambda \in \Z_+ \}$. More explicitly, the
set of Pareto-optimal solutions for $MIP_{A,C}(b)$ is the set of
minimal elements (with respect to $\prec_C$) of $\Gamma$.
\end{remark}

In order to illustrate the above algorithm, we present an example
of MOILP with two objectives where all the computations are done
in detail.

\begin{ej}
\label{ejfinal}

\begin{equation}
\label{multiobj}
\begin{array}{lrl}
\min & \;  \{10x + y , x + 10y\}& \\
s.a.& &\\
 & 2x + 2y & \geqslant 17\\
 & 2y & \leqslant 11\\
 & x & \leqslant 10\\
 & x, y &\in \mathbb{Z}_+
 \end{array}
 \end{equation}

Transforming the problem to the standard form results in:

\begin{equation}
\label{multiobj2}
\begin{array}{lrl}
\min & \;  \{10x + y + 0z + 0t + 0q, x + 10y + 0z + 0t + 0q\}& \\
s.a.& &\\
 & 2x + 2y - z & = 17\\
 & 2y + t &= 11\\
 & x + q &= 10\\
 & x, y, z, t, q&\in \mathbb{Z}_+
 \end{array}
 \end{equation}

\begin{description}
\item[\textsc{Step 1. }] Feasible solution for $MIP_{A,C}(b)$: $u = (9,4,9,3,1)$.
\item[\textsc{Step 2. }] Following the steps of Algorithm \ref{systemofgenerators}:
 \begin{enumerate}
\item Basis for $Ker(A)$ : $\mathcal{B} := \{(0,1,2,-2,0),(-1,0,-2,0,1)\}$.
\item LLL-reduced basis for $\mathcal{B}$ : $\mathcal{B}_{red} := \mathcal{B} :=
\{(-1,0,-2,0,1),(-1,1,0,-2,1)\}$.
\item $J_0 := \langle x^{u_+} - x^{u_-} : u \in \mathcal{B}_{red}
\rangle  =\langle x_5-x_1x_3^2, x_2x_5 - x_1x_4^2 \rangle$
\item $J_{i+1} := (J_{i}:x_i^\infty)$
\begin{enumerate}
    \item $\widetilde{\grob}_0 := \{x_5-x_1x_3^2, x_2x_5 - x_1x_4^2, x_2x_3^2-x_4^2\} \Rightarrow
    J_1 :=\langle  x_5-x_1x_3^2, x_2x_5 - x_1x_4^2, x_2x_3^2-x_4^2\rangle$
    \item $\widetilde{\grob}_1 := \{x_5-x_1x_3^2, x_2x_5 - x_1x_4^2, x_2x_3^2-x_4^2\} \Rightarrow
    J_2 :=\langle  x_5-x_1x_3^2, x_2x_5 - x_1x_4^2, x_2x_3^2-x_4^2\rangle$
        \item $\widetilde{\grob}_2 := \{x_5-x_1x_3^2, x_2x_5 - x_1x_4^2, x_2x_3^2-x_4^2\} \Rightarrow
    J_3 :=\langle  x_5-x_1x_3^2, x_2x_5 - x_1x_4^2, x_2x_3^2-x_4^2\rangle$
        \item $\widetilde{\grob}_3 := \{x_5-x_1x_3^2, x_2x_5 - x_1x_4^2, x_2x_3^2-x_4^2\} \Rightarrow
    J_4 :=\langle  x_5-x_1x_3^2, x_2x_5 - x_1x_4^2, x_2x_3^2-x_4^2\rangle$
\end{enumerate}
\item $\Im_A = \langle x_5-x_1x_3^2, x_2x_5 - x_1x_4^2, x_2x_3^2-x_4^2,
x_1x_3^2-1\rangle \mapsto$ \\
$I_A = \langle \{ \big( (1,0,0,0,1),(0,1,0,2,0)\big), \big(
(1,0,2,0,0),(0,0,0,0,1)\big), \big( (0,1,2,0,0),
(0,0,0,2,0)\big)\}$
\end{enumerate}

\item[\textsc{Step 3. }] Computing a p-Gröbner basis for $I_A$, using the order
$\prec^s_C$ (Remark \ref{remark:slacks}), and following Algorithm
\ref{pbuchberger2} we obtain $\grob$, whose maximal chains are:

\begin{description}
    \item[$\grob_1$] $\{ \big( (0,1,2,0,0), (0,0,0,2,0), (0,1,2,0,0)\big), \big( (0,1,0,0,2),(2,0,2,2,0),(0,1,0,0,2)\big),$\\
    $\big( (0,1,0,0,1),(1,0,0,2,0),(0,1,0,0,1)\big)\}$.
    \item[$\grob_2$] $\{ \big( (1,0,0,4,0), (0,2,2,0,1), (1,0,0,4,0) \big), \big( (1,0,2,0,0), (0,0,0,0,1), (1,0,2,0,0)\big),$\\
    $\big( (1,0,0,2,0), (0,1,0,0,1), (1,0,0,2,0) \big)\}$.
\end{description}
\item[\textsc{Step 4. }] Partial remainders:
Reducing first by $\grob_1$:
\begin{itemize}
\item[] $pRem((9,4,9,3,1), \grob_1) = \{(9,0,1,11,1)\}$.\\
Then, reducing each remainder by $\grob_2$:
\item[] $pRem((9,0,1,11,1), \grob_2) =   \{(9,0,1,11,1), (8,2,3,7,2),
(7,2,1,9,3),
(6,3,1,5,4), (5,4,1,3,5), (4,5,1,1,6)\}$.\\
\end{itemize}
The entire set of Pareto-optimal solutions is:
$$
\{(9,0,1,11,1),(8,1,1,9,2),(7,2,1,7,3),(6,3,1,5,4),(5,4,1,3,5),(4,5,1,1,6)\}
$$
\end{description}
Figure \ref{example1-graph} shows the feasible region and the
Pareto-optimal solutions of the example above.
\begin{figure}[h]

\begin{center}
\begin{pspicture}(0,0)(10.,6.)
\newgray{grismedio}{0.8}
\pspolygon*[linecolor=yellow](3,5.5)(10,5.5)(10,0)(8.5,0)
  \psgrid[subgriddiv=0,griddots=10,gridlabels=7pt]
  \psaxes[linewidth=1pt,%
    ticks=none,%
    labels=none]{->}(0,0)(0,0)(10.5,6.5)
  \uput[-30](10,-0.5){$\mathbf{x}$}
  \uput[120](-0.5,6){$\mathbf{y}$}

    \psline[linewidth=1pt, linecolor=red](3,5.5)(10,5.5)
    \psline[linewidth=1pt, linecolor=red](3,5.5)(8.5,0)
    \psline[linewidth=1pt, linecolor=red](10,0)(8.5,0)
    \psline[linewidth=1pt, linecolor=red](10,5.5)(10,0)
    \psdots[dotsize=4pt, linecolor=blue](9,0)(8,1)(7,2)(6,3)(5,4)(4,5)
    \pscircle[linewidth=0.2pt](10,0){0.08}
        \pscircle[linewidth=0.2pt](10,0){0.08}
        \pscircle[linewidth=0.2pt](10,1){0.08}
        \pscircle[linewidth=0.2pt](10,2){0.08}
        \pscircle[linewidth=0.2pt](10,3){0.08}
        \pscircle[linewidth=0.2pt](10,4){0.08}
        \pscircle[linewidth=0.2pt](10,5){0.08}
        \pscircle[linewidth=0.2pt](9,1){0.08}
        \pscircle[linewidth=0.2pt](9,2){0.08}
        \pscircle[linewidth=0.2pt](9,3){0.08}
        \pscircle[linewidth=0.2pt](9,4){0.08}
        \pscircle[linewidth=0.2pt](9,5){0.08}
        \pscircle[linewidth=0.2pt](8,2){0.08}
        \pscircle[linewidth=0.2pt](8,3){0.08}
        \pscircle[linewidth=0.2pt](8,4){0.08}
        \pscircle[linewidth=0.2pt](8,5){0.08}
        \pscircle[linewidth=0.2pt](7,3){0.08}
        \pscircle[linewidth=0.2pt](7,4){0.08}
        \pscircle[linewidth=0.2pt](7,5){0.08}
        \pscircle[linewidth=0.2pt](6,4){0.08}
        \pscircle[linewidth=0.2pt](6,5){0.08}
        \pscircle[linewidth=0.2pt](5,5){0.08}

    \psline[linecolor=blue]{->}(2,2)(1,2.1)
    \psline[linecolor=blue]{->}(2,2)(2.1,1)
\end{pspicture}
\vspace*{0.5cm}
\end{center}
\caption{\label{example1-graph}Feasible region, Pareto-optimal
solutions and improvement cone for Example \ref{ejfinal}}
\end{figure}
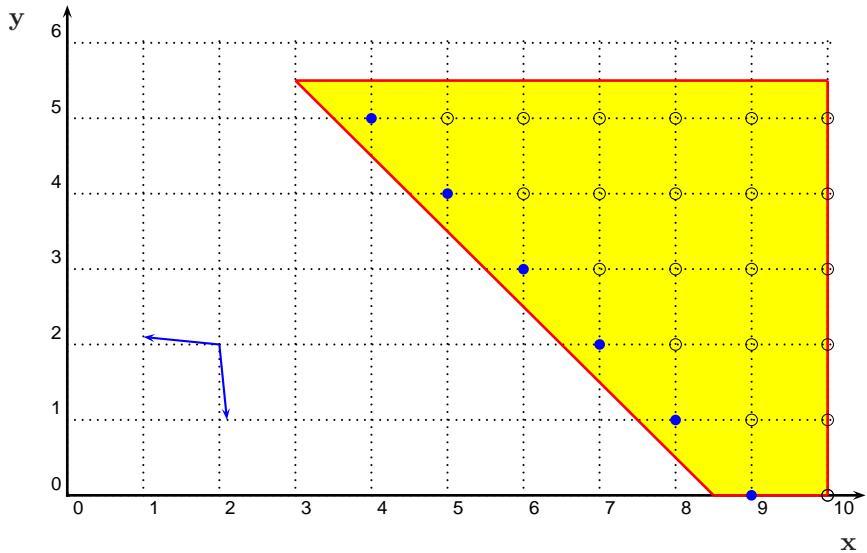

In addition, we have evaluated the problem with the same feasible
region but choosing a cost matrix such that the respective normal
vectors of each of the rows in the matrix form an acute angle.
Then, non
supported solutions appear in the set of Pareto-optimal solutions. Figure \ref{example2-graph} shows the Pareto-optimal solutions for the same feasible region and $C=\left[ \begin{array}{cc} 10 & -1 \\
-1 & 10
\end{array}\right]$.

\begin{figure}[h]
\vspace*{0.3cm}
\begin{center}
\begin{pspicture}(0,0)(10.,6.)
\newgray{grismedio}{0.8}
\pspolygon*[linecolor=yellow](3,5.5)(10,5.5)(10,0)(8.5,0)
  \psgrid[subgriddiv=0,griddots=10,gridlabels=7pt]
  \psaxes[linewidth=1pt,%
    ticks=none,%
    labels=none]{->}(0,0)(0,0)(10.5,6.5)
  \uput[-30](10,-0.5){$\mathbf{x}$}
  \uput[120](-0.5,6){$\mathbf{y}$}

    \psline[linewidth=1pt, linecolor=red](3,5.5)(10,5.5)
    \psline[linewidth=1pt, linecolor=red](3,5.5)(8.5,0)
    \psline[linewidth=1pt, linecolor=red](10,0)(8.5,0)
    \psline[linewidth=1pt, linecolor=red](10,5.5)(10,0)
    \psdots[dotsize=4pt,
    linecolor=blue](9,1)(8,1)(7,2)(6,3)(5,4)(4,5)(5,5)(6,4)(7,3)(8,2)(9,0)
    (10,0)
    \pscircle[linewidth=0.2pt](10,0){0.08}
        \pscircle[linewidth=0.2pt](10,1){0.08}
        \pscircle[linewidth=0.2pt](10,2){0.08}
        \pscircle[linewidth=0.2pt](10,3){0.08}
        \pscircle[linewidth=0.2pt](10,4){0.08}
        \pscircle[linewidth=0.2pt](10,5){0.08}
        \pscircle[linewidth=0.2pt](9,2){0.08}
        \pscircle[linewidth=0.2pt](9,3){0.08}
        \pscircle[linewidth=0.2pt](9,4){0.08}
        \pscircle[linewidth=0.2pt](9,5){0.08}
        \pscircle[linewidth=0.2pt](8,3){0.08}
        \pscircle[linewidth=0.2pt](8,4){0.08}
        \pscircle[linewidth=0.2pt](8,5){0.08}
        \pscircle[linewidth=0.2pt](7,4){0.08}
        \pscircle[linewidth=0.2pt](7,5){0.08}
        \pscircle[linewidth=0.2pt](6,5){0.08}

    \psline[linecolor=blue]{->}(2,2)(1,1.9)
    \psline[linecolor=blue]{->}(2,2)(1.9,1)
\end{pspicture}
\vspace*{0.5cm}
\end{center}
\caption{\label{example2-graph}Feasible region, Pareto-optimal
solutions and improvement cone for Example \ref{ejfinal} with
$C=\left[ [10, -1],[-1, 10]\right]$}
\end{figure}
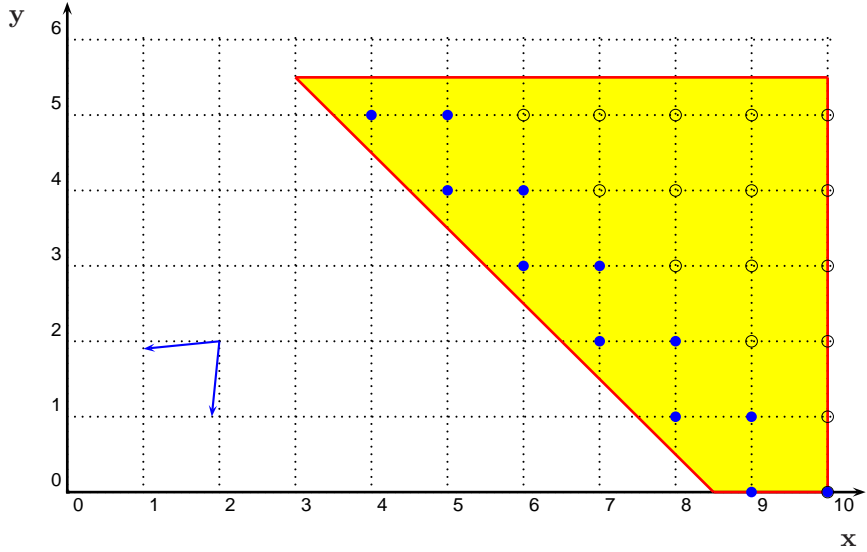

\end{ej}
\setcounter{section}{3}

\section{Computational Results}
\label{sect3} A series of computational experiments have been
performed in order to evaluate the behavior of the proposed
solution method. Programs have been coded in MAPLE 10 and executed
in a PC with an Intel Pentium 4 processor at 2.66GHz and 1 GB of
RAM. In the implementation of Algorithm \ref{pbuchberger2} to
obtain the p-Gr\"obner basis, the package \textit{poset} for Maple
\cite{stembridge06} has been used to compute, at each iteration,
the maximal chains for the p-Gr\"obner basis. The implementation
has been done in a symbolic programming language, available upon
request, in order to make the access easy to both optimizers and
algebraic geometers.

The performance of the algorithm was tested on randomly generated
instances for knapsack and transportation multiobjective problems
for 2, 3 and 4 objectives. For the knapsack problems, 4, 5 and 6
variables programs were considered, and for each group, the
coefficients of the constraint were randomly generated in $[0,20]$
and the coefficients of the objective matrices range in $[0,20]$.
Once the constraint vector, $(a_1, \ldots, a_n)$, is generated,
the right-hand side is fixed as $b=\lceil \frac{1}{2} \sum_{i=1}^n
a_i \rceil$ to ensure feasibility.

The computational tests have been done in the following way for
each number of variables: (1) Generate 5 constraint vectors and
compute the initial system of generators for each of them using
Algorithm \ref{systemofgenerators}; (2) Generate five random
objective matrices for each number of objectives (2, 3 and 4) and
compute the corresponding p-Gr\"obner basis using Algorithm
\ref{pbuchberger2}; and (3) with $b=\lceil \frac{1}{2}
\sum_{i=1}^n a_i \rceil$ and for each objective matrix, compute
the Pareto-optimal solutions using Algorithm \ref{pos}.

Table \ref{comp:knapsack} contains a summary of the average
results obtained for the considered knapsack multiobjective
problems. The second, third and fourth columns show the average
CPU times for each stage in the algorithm: \verb"sog" is the CPU
time for computing the system of generators, \texttt{pgr\"obner}
is the CPU time for computing a p-Gr\"obner basis, and
\texttt{pos} is the time for computing a feasible solution and
partially reduce it to obtain the set of Pareto-optimal solutions.
The fifth column shows the total time for computing the set of
Pareto-optimal solutions for the problem. Finally, the sixth and
seventh columns show the average number of Pareto-optimal
solutions and the number of maximal chains in the p-Gr\"obner
basis for the problem. The problems have been named as
\verb"knapN_O" where \verb"N" is the number of variables and
\verb"O" is the number of objectives.

\begin{table}[h]\small{\center{
\begin{tabular}{lrrrrrrrr}\hline
Problem & \texttt{sog} & \texttt{pgr\"obner} &
\texttt{pos}& total & $|POS|$ & $|MaxChains|$& \verb"steps" & \verb"act_pGB" \\
\hline

\verb"knap4_2" & 0.063 & 249.369 & 1.265 & 250.697 & 11& 20 & 2 & 164.920\\
\verb"knap4_3" & 0.063 & 1002.689 & 2.012  & 1004.
704 & 5 & 46& 2 & 772.772\\
\verb"knap4_4" & 0.063 & 1148.574 & 2.374 & 1151.011 & 16 & 98& 2.4 & 763.686\\
\verb"knap5_2" & 0.125 & 1608.892 & 0.875 & 609.892 & 3 & 29& 2 & 1187.201\\
\verb"knap5_3" & 0.125 & 3500.831 & 2.035 & 3503.963 & 2 & 30& 2.2 & 2204.123\\
\verb"knap5_4" & 0.125 & 3956.534 & 2.114 & 3958.773 & 9 & 45.4& 3 & 3044.157\\
\verb"knap6_2" & 0.185 & 2780.856 & 2.124 & 2783.165 & 18& 156& 2.4 & 2241.091\\
\verb"knap6_3" & 0.185 & 3869.156 & 2.018 & 3871.359 & 16.4& 189& 2.4 & 2790.822\\
\verb"knap6_4" & 0.185 & 4598.258 & 3.006 & 4601.449 & 26& 298& 3.2 & 3096.466\\
\hline\\
&&&&&&
\end{tabular} \caption{\label{comp:knapsack}Summary of computational
experiments for knapsack problems}}}
\end{table}

For the transportation problems, instances with 3 origins $\times$
2 destinations, 3 origins $\times$ 3 destinations and 4 origins
$\times$ 2 destinations were considered. In this case, for fixed
numbers of origins, $s$, and destinations, $d$, the constraint
matrix, $A \in Z^{(s+d)\times (sd)}$, is fixed. Then, we have
generated 5 instances for each problem of size $s\times d$. Each
of these instances is combined with 5 different right-hand side
vectors. The procedure is analogous to the knapsack computational
test: a first step where a system of generators is computed, a
second one, where the p-Gr\"obner basis is built and in the last
step, the set of Pareto-optimal solutions is computed using
partial reductions. Table \ref{comp:trans} shows the average CPU
times and the average number of Pareto-optimal solutions and
maximal chains in the p-Gr\"obner basis for each problem. The
\verb"step" column shows the average number of steps in the
p-Gr\"obner computation, and \verb"act_pGB" is the average CPU
time in the computation of the p-Gr\"obner basis elapsed since the
last element was added to the basis until the end of the process.
The problems have been named as \verb"transNxM_O" where \verb"N"
is the number of origins, \verb"M" is the number of destinations
and \verb"O" is the number of objectives.

\begin{table}[h]\small{\center{
\begin{tabular}{lrrrrrrrr}\hline
Problem & \texttt{sog} & \texttt{pgr\"obner} &
\texttt{pos}& total & $|POS|$ & $|MaxChains|$ & \verb"steps" & \verb"act_pGB"\\
\hline
\verb"tranp3x2_2" & 0.015 & 11.813 & 0.000 & 11.828& 5.2 & 6 & 2 & 7.547\\
\verb"tranp3x2_3" & 0.015 & 7.218 & 13.108 & 30.341 & 12 & 2.6& 2 & 6.207\\
\verb"tranp3x2_4" & 0.015 & 6.708 & 15.791 & 21.931 & 6 & 5& 2.2 & 4.561\\
\verb"tranp3x3_2" & 0.047 & 1545.916 & 1.718 & 1547.681 & 5 & 92& 2 & 928.222\\
\verb"tranp3x3_3" & 0.047 & 3194.333 & 11.235 & 3205.615 & 9& 122& 2.4 & 2172.146\\
\verb"tranp3x3_4" & 0.047 & 3724.657 & 7.823 & 3732.527 & 24 & 187.4& 2.2 & 2112.287\\
\verb"tranp4x2_2" & 0.046 & 675.138 & 2.122 & 677.306 & 3.4 & 35.2& 2 & 398.093\\
\verb"tranp4x2_3" & 0.046 & 1499.294 &6.288 & 1505.628 & 5.8 & 42.4 & 2.2 & 119.519\\
\verb"tranp4x2_4" & 0.046 & 2285.365 &7.025 & 2292.436 & 12 & 59 & 2.2 & 1654.048\\
\hline\\
&&&&&&
\end{tabular}

\caption{\label{comp:trans}Summary of computational experiments
for the battery of multiobjective transportation problems}}}

\end{table}

As can be seen in tables \ref{comp:knapsack} and \ref{comp:trans},
the overall CPU times are clearly divided into the three steps,
being the most costly the computation of the p-Gr\"obner basis. In
all the cases more than 99\% of the total time is spent computing
the p-Gr\"obner basis. Once this structure is computed, obtaining
the Pareto-optimal solutions is done very efficiently.

The CPU times and sizes in the different steps of the algorithm
are highly sensitive to the number of variables. However, our
algorithm is not very sensitive to the number of objectives, since
the increment of CPU times with respect to the number of
objectives is much smaller than the one with  respect to the
number of variables.

It is clear that one can not expect fast algorithms for solving
MOILP, since all these problems are NP-hard. Nevertheless, our
approach gives exact tools that apart from solving these problems,
give insights into the geometric and algebraic nature of the
problem.

As mention above, using our methodology one can identify the
common algebraic structure within any multiobjective integer
linear problem. This connection allows to improve the efficiency
of our algorithm making use of any advance that improves the
computation of Gr\"obner bases. In fact, any improvements of the
standard Gr\"obner bases theory may have an impact in improving
the performance of this algorithm. In particular, one can expect
improvements in the efficiency of our algorithm based on the
special structure of the integer program (see for instance Remark
\ref{remark:dimkernel}). In addition, we have to mention another
important issue in our methodology. As shown in Theorem
\ref{testfam-pgrobner}, solving MOILP with the same constraint and
objective matrices requires computing only once the p-Gr\"obner
basis. Therefore, once this is done, we can solve different
instances varying the right-hand side very quickly.

Finally, we have observed from our computational tests that a
significant amount of the time, more than 60\% of the time (see
column \verb"act_pGB"), for the computation of the p-Gr\"obner
basis is spent checking that no new elements are needed in this
structure. This implies that the actual p-Gr\"obner basis is
obtained much earlier than when the final test is finished. A
different truncation strategy may be based on the number of steps
required to obtain the p-Gr\"obner basis. According to the exact
method, the algorithm stops once in a step no new elements are
added to the structure. Our tables show that in most cases the
number of steps is $2$, actually only one step is required to
generate the entire p-Gr\"obner basis (see column \verb"steps").
These facts can be used to accelerate the computational times at
the price of obtaining only heuristic Pareto-optimal solutions.
This idea may be considered an alternative primal heuristic in
MOILP and will be the subject of further research.


\begin{thebibliography}{99}

\bibitem{aardal02}  Aardal, K. , Weismantel, R., Wolsey,
L. (2002), Non-Standard approaches to integer programming,
Discrete Applied Mathematics 123  5--74.

\bibitem{adams94} Adams, W. , Loustaunau, P. (1994), An
introduction to Gr\"obner bases. Graduate Studies in Mathematics
3. American Mathematical Society.

\bibitem{bazaraa93}  Bazaraa, M.S, Sherali, H.D. and Shetty, C.M. (1993). Nonlinear programming : theory and
algorithms. New York [etc.]: John Wiley and Sons.

\bibitem{baer69}  Baer, R.M,  and \O sterby, O. (1969), Algorithms over Partially Ordered
Sets. Journal BIT Numerical Mathematics vol 9, 2, 97--118.


\bibitem{bertsimas05} Bertsimas, D. and Weismantel, R. (2005).
Optimization over integers. Dynamic Ideas, Belmont, Massachusetts.
ISBN: 0-9759146-2-6.

\bibitem{buchberger65} Buchberger, B. (1965), An Algorithm for Finding a Basis for the Residue Class Ring of a Zero-Dimensional Polynomial
Ideal.PhD thesis, University of Innsbruck, Institute for
Mathematics.

\bibitem{cantor1897} Cantor, G. (1897), Beitr\"{a}ge zur Begr\"{u}ndung der transfiniten Mengenlehre (Zweiter Artikel), \emph{Math.\ Ann.}\ 49, 207--246.

\bibitem{cayley1889} Cayley, A. (1889), A theorem on trees. Quarterly Journal of Mathematics 23, 376-378.
\bibitem{chankong-haimes83} Chankong, V. and Haimes,Y.Y.  (1983). Multiobjective Decision Making Theory and Methodology. Elsevier Science, New
York.
\bibitem{conti-traverso91} Conti, P. and Traverso, C. (1991), Buchberger algorithm and
integer programming, Proceedings AAECC9 (New Orleans), Springer
LNCS 539, 130 - 139.

\bibitem{cox-little-oshea92} Cox, D., Little, J., OShea, D. (1992), Ideals, varieties, and algorithms : an introduction to Computational Algebraic Geometry and Commutative Algebra. first edition. Springer-Verlag, New York.

\bibitem{cox-little-oshea98} Cox, D., Little, J., OShea, D. (1998), Using Algebraic Geometry. first edition. Springer-Verlag, New York.

\bibitem{delorme03}Delorme, X. , Gandibleux, X. and Degoutin, F. (2003).
Resolution approch\'e du probleme de set packing bi-objectifs. In
Proceedings de l'ecole d'Automne de Recherche Operationnelle de
Tours (EARO), 74--80.

\bibitem{edgeworth1881} Edgeworth, F.Y. (1881), Mathematical Psychics. P. Keagan, London.

\bibitem{ehrgott00} Ehrgott, M. (2000). Multicriteria optimization. Lecture Notes in Economics
and Mathematical Systems 491. Springer-Verlag, Berlin.

\bibitem{ehrgott-gandibleux00}Ehrgott, M. and Gandibleux, X. (2000) A survey and annotated bibliography
of multicriteria combinatorial optimization. OR Spektrum
22:425-460.

\bibitem{ehrgott02} Ehrgott, M. and Gandibleux, X. (editors) (2002). Multiple Criteria
Optimization. State of the Art Annotated Bibliographic Surveys.
Boston, Kluwer.

\bibitem{ehrgott-ryan02}Ehrgott, M. and Ryan, D.M. (2002). Constructing robust crew
schedules with bicriteria optimization. Journal of Multi-Criteria
Decision Analysis 11(3):139--150.

\bibitem{ehrgott05} Ehrgott, M., Figueira, J. and Greco, S. (editors) (2005). Multiple Criteria
Decision Analysis. State of the Art Surveys. New York, Springer.

\bibitem{ehrgott06} Ehrgott, M. ,  Figueira, J. and Gandibleux, X. (editors) (2006). Multiobjective
Discrete and Combinatorial Optimization. Annals of Operations
Research 147.

\bibitem{faugere99} Faug\`{e}re, J.C. (1999). A new efficient algorithm for computing Grobner bases
(F4). Journal of Pure and Applied Algebra, Volume 139, Number 1,
June 1999 , pp. 61--88(28).

\bibitem{faugere02} Faug\`{e}re, J.C. (2002). A new efficient algorithm for computing Gröbner bases without reduction to zero F5. In T. Mora, editor, Proceedings of the 2002 International Symposium on Symbolic and Algebraic Computation ISSAC,
pages 75--83. ACM Press, July 2002. ISBN: 1-58113-484-3.

\bibitem{fernandez-puerto03}  Fern\'andez, E. and Puerto, J. (2003). The multiobjective solution of the uncapacitated plant location problem''. European Journal of Operational Research. vol. 45, n.3 509-529.

\bibitem{hamacher94} Hamacher, H. and Ruhe, G. (1994).  On spanning tree problems with multiple objectives. Annals of Operations
Research 52(4) 209--230.

\bibitem{hausdorff1906} Hausdorff, F. (1906),
Untersuchungen über Ordungtypen, Berichte über die Verhandlungen
der königlich sächsischen Gesellschaft der Wissenschaften zu
Leipzig, Matematisch - Physische Klasse 58, pp. 106--169.

\bibitem{hosten-sturmfels95} Hosten, S. and Sturmfels, B. (1995) GRIN: An implementation of Gr\"obner bases for integer programming, Integer Programming and Combinatorial Optimization (E. Balas and J. Clausen, eds.), Lecture Notes in Computer Science, no. 920, Springer, Berlin, 1995, Proceedings of the 4th International IPCO Conference, pp. 267 -- 276.

\bibitem{hosten97} Hosten, S. (1997)  Degrees of Gr\"obner bases of integer
programs, Ph.D Thesis Cornell University.

\bibitem{hosten-thomas98} Hosten, S. and Thomas, R. R. (1998), Gr\"obner bases and integer programming, In:
Gr\"obner Bases and Applications, B. Buchberger, F. Winkler
(Eds.), London Mathematical Society Lecture Note Series 251,
Cambridge University Press, Cambridge UK, 1998, 470--484.

\bibitem{ishibuchi-murata98} Ishibuchi, H. and Murata, T. (1998),  A multi-objective genetic local search algorithm and its application to flowshop scheduling, IEEE Trans. Syst., Man, Cybern. C 28, 392--403.

\bibitem{jozefowiez04} Jozefowiez, N. , Semet, F. and Talbi, E-G. (2004). A multi-objective evolutionary algorithm for the covering tour problem, Chapter 11 in "Applications
 of multi-objective evolutionary algorithms", C. A. Coello and G. B. Lamont (editors), p 247-267, World Scientific.

\bibitem{karwan-villareal82} Karwan, M.H. and Villarreal, B. (1982).
Multicriteria Dynamic Programming with an Application to the
Integer Case. Journal of Optimization Theory and Applications.
Vol. 31. pp 43-69.


\bibitem{lll82} Lenstra, A. K. , Lenstra, H. W. and Lov\`asz, L. (1982),
Factoring Polynomials with Rational Coefficients, Math. Ann. 261,
4.


\bibitem{miettinen99} Miettinen, K. (1999). Nonlinear Multiobjective Optimization,
Kluwer Academic Publishers, Boston. ISBN: 0792382781.

%
\bibitem{pareto1896} Pareto, V. (1896), Manual d'economie politique, F. Rouge, Lausanne.

\bibitem{pottier91} Pottier, L. (1991), Minimal solutions of linear diophantine systems: bounds and algorithms. In Proc. of the Fourth Intern. Conf. on Rewriting Techniques and Applications
162-173.

\bibitem{tanino85} Sawaragi, Y., Nakayama, H. and Tanino, T. (1985). Theory of Multiobjective Optimization, Academic
Press.

\bibitem{schrijver03} Schrijver, A. (2003), Combinatorial Optimization: Polyhedra and Efficiency. New York, NY: Springer-Verlag, 2003. ISBN: 3540443894.

\bibitem{sherbeny01} El-Sherbeny, N. (2001). Resolution of a Vehicle Routing
Problem with Multiobjective Simulated Annealing Method, PhD
thesis, Faculte Polytechnique de Mons, Belgium.

\bibitem{sourd06} Sourd, F.,  Spanjaard, O. and Perny, P. (2006). Multi-objective branch and bound. Application to the bi-objective spanning tree problem. In 7th International Conference in Multi-Objective Programming and Goal
Programming. Loire Valley (City of Tours), France.

\bibitem{stembridge06} Stembridge, J.R. (2006). A Maple Package
for Posets. Available at \verb"www.math.lsa.umich.edu/~jrs".

\bibitem{steuer85} Steuer, R.E. (1985). Multiple Criteria Optimization: Theory,
Computation and Application. John Wiley \& Sons, New York, NY.

\bibitem{sturmfels95} Sturmfels, B. (1996), Gr\"obner Bases and Convex Polytopes, American Mathematical Society,
University Lectures Series, No. 8, Providence, Rhode Island.

\bibitem{sturmfels02} Sturmfels, B. (2002), Solving Systems of Polynomial Equations. CBMS Regional Conference Series of the AMS, Number
97, Rhode Island.

\bibitem{sturmfels04} Sturmfels, B. (2004), Algebraic Recipes for integer
programming. Proceedings of Symposia in Applied Mathematics. ISSN
: 0160-7634. Vol. 61

\bibitem{sturmfels-thomas97} Sturmfels, B. and Thomas, R. R.
(1997), Variation of Cost Functions in Integer Programming.
Mathematical Programming 77, 357--387.

\bibitem{thomas95} Thomas, R. R. (1995), A geometric Buchberger algorithm for integer programming, Mathematics of Operations Research 20, 4 (1995).

\bibitem{thomas97} Thomas, R. R. (1997), Applications to Integer Programming PSAM: Proceedings of the 53th Symposium in Applied
Mathematics, American Mathematical Society, 119-142.



\bibitem{thomas-weismantel97} Thomas, R. R., Weismantel, R. (1997), Truncated Groebner bases for integer programming, Applicable Algebra in Engineering, Communication and Computing 8, 241-257.


\bibitem{urbaniak97} Urbaniak, R., Weismantel, R., Ziegler, G.M. (1997)  A variant of the Buchberger algorithm for integer programming, SIAM J. Discrete Math 10 96-108.
\bibitem{yu74} Yu, P.L., Cone Convexity, Cone Extreme Points and Nondominated Solutions in Decision Problems with Multiobjectives. Journal of Optimization Theory and Applications 14 \#3 319-377.

\bibitem{zionts79} Zionts, S. (1979). A survey of multiple criteria integer programming methods. Annals of Discrete Mathematics 5, 389--398.
\bibitem{zionts-wallenius80} Zionts, S., Wallenius, J. (1980). Identifying efficient vectors:
some theory and computational results. Operations Research 23,
785--793.

\end{thebibliography}
\end{document}